\input amstex
\documentstyle{amsppt}
\magnification=1200
\overfullrule=0pt
\loadbold

\TagsOnLeft
\TagsAsMath

\redefine\leq{\leqslant}
\redefine\geq{\geqslant}

\define\VV{{\Cal V}}
\define\WW{{\Cal W}}
\define\UU{{\Cal U}}
\define\MM{{\Cal M}}
\define\HH{{\Cal H}}
\define\GG{{\Cal G}}

\define\gt{{\Cal G}t}

\define\aand{\quad\text{and}\quad}
\define\wwhere{\quad\text{where}\quad}
\define\ffor{\quad\text{for}\quad}
\define\for{~\text{for}~}

\define\iiif{~\text{if}~}

\define\TTT{{\bold T}}
\define\RRR{{\bold R}}
\define\ZZZ{{\bold Z}}
\define\NNN{{\bold N}}
\define\KKK{{\bold K}}

\define\kk{{\Bbb K}}
\define\Wkr{W^{\circ}}

\define\ind{\text{\rm ind}}

\redefine\Im{\text{\rm Im }}
\define\supp{\text{\rm supp }}
\define\Int{\text{\rm Int }}

\define\nr{\Vert}
\define\smo{C^{\infty}}

\define\fpr#1#2{{#1}^{-1}({#2})}
\define\sdvg#1#2#3{\widehat{#1}_{[{#2},{#3}]}}
\define\disc#1#2#3{B^{({#1})}_{#2}({#3})}
\define\Disc#1#2#3{D^{({#1})}_{#2}({#3})}
\define\desc#1#2#3{B_{#1}(\leq{#2},{#3})}
\define\Desc#1#2#3{D_{#1}(\leq{#2},{#3})}
\define\komp#1#2#3{{\bold
K}({#1})^{({#2})}({#3})}
\define\Komp#1#2#3{\big({\bold
K}({#1})\big)^{({#2})}({#3})}

\define\ran{\{(A_\lambda , B_\lambda)\}_{\lambda\in\Lambda}}
\define\rans{\{(A_\sigma , B_\sigma)\}_{\sigma\in\Sigma}}

\define\fmin{F^{-1}}

\define\chart{\Phi_p:U_p\to B^n(0,r_p)}
\define\atlas{\{\Phi_p:U_p\to B^n(0,r_p)\}_{p\in S(f)}}
\define\flow{{\Cal V}=(f,v,{\Cal U})}

\define\LLN{{\Cal {LN}}}

\define\capp{\pitchfork}

\define\crr{p\in S(f)}
\define\nrv{\Vert v \Vert}
\define\nrw{\Vert w \Vert}

\define\obb{\cup_{p\in S(f)} U_p}

\define\stind#1#2#3{{#1}^{\displaystyle\rightsquigarrow}_
{[{#2},{#3}]}}

\define\indl#1{{\scriptstyle{\text{\rm ind}\leqslant {#1}~}}}
\define\inde#1{{\scriptstyle{\text{\rm ind}      =   {#1}~}}}

\define\xx{\Bbb X}
\define\aaaa{\Bbb A}

\define\yy{\Bbb Y}
\define\tttt{\Bbb T}

\define\XXX{\bold X}

\define\AAA{\bold A}

\define\xxx{\Bbb X =\{X_0,...,X_k\}}

\define\vem{\text{Vectt}(M)}

\define\id{\text{id}}

\define\aaa{\Bbb A}

\define\st#1{\overset\rightsquigarrow\to{#1}}
\define\bst#1{\overset{\displaystyle\rightsquigarrow}\to{\boldkey{#1}}}

\define\stexp#1{{#1}^{\rightsquigarrow}}
\define\bstexp#1{{#1}^{\displaystyle\rightsquigarrow}}

\define\bstind#1#2#3{{\boldkey{#1}}^{\displaystyle\rightsquigarrow}_
{[{#2},{#3}]}}
\define\bminstind#1#2#3{\stind{({\boldkey{-}\boldkey{#1}})}{#2}{#3}}

\define\Tb{\text{ \rm Tb}}
\define\llll{\Bbb L}

\define\VODIN{V_{1/3}}
\define\VDVA{V_{2/3}}
\define\VM{V_{1/2}}

\define\where{\quad\text{\rm where}\quad}

\define\kr#1{{#1}^{\circ}}

\define\TT{\Cal T}

\define\bb{{\Bbb B}}

\define\vew{\text{\rm Vect}^1 (W,\bot)}

\define\BBB{\bold B}

\define\hrrr{\text{\rm Vect}^1(M)}
\define\vemm{\text{\rm Vect}^1_0(M)}

\define\verr{\text{\rm Vect}^1_0(\RRR^ n)}

\define\mods{\vert s(t)\vert}

\topmatter
\title
The incidence coefficients
in the Novikov Complex are generically rational
functions
\endtitle
\leftheadtext{A.V.Pajitnov}
\rightheadtext{      Incidence coefficients in the
Novikov Complex}
\author
A.V.Pajitnov
\endauthor
\address
Universit\'e de Nantes, Facult\'e de sciences, 2, rue de la
Houssini\`ere, 44072, Nantes, Cedex
\endaddress
\email
pajitnov\@univ-rennes1.fr
\endemail

\abstract
For a Morse map $f:M\to S^1$ Novikov [11]
has introduced an analog of  Morse complex,
defined over the ring $\ZZZ[[t]][t^{-1}]$
of  integer Laurent power series.  Novikov conjectured, that
generically the matrix entries of the differentials
in this complex are of the form
$\sum_ia_it^i$, where
$a_i$
grow at most exponentially in $i$.
We prove that for any given $f$ for a $C^0$ generic
gradient-like vector field all the incidence coefficients
above are rational functions in $t$
(which implies obviously the exponential growth rate estimate).
\endabstract
\endtopmatter

\head{ Introduction} \endhead
\subhead{A.  Morse-Novikov theory}\endsubhead
 The classical Morse-Thom-Smale construction
associates to a Morse function $f:M\to\RRR$
on a closed manifold a free chain complex
$C_*(f)$ where the number $m(C_p(f))$ of free generators
of $C_p(f)$ equals the number of the critical points of $f$ of
index $p$ for each $p$. The boundary operator in this
complex is defined in a geometric way, using the trajectories
of a gradient of $f$, joining critical points of
$f$ (see [8], [21], [23], [25]).

In the early 80s S.P.Novikov generalized this construction
 to the case of maps $f:M\to S^1$ (he was lead
to this generalization by considering a problem
of motion of a solid in a fluid, see [9], [10]).
The corresponding analog of Morse complex
is a free chain complex $C_*(f)$ over
$\ZZZ [[t]][t^{-1}]$. Its number of free generators
equals the number of critical points of $f$ of index $p$,
and the homology of $C_*(f)$ equals to the completed
homology of the cyclic covering. In [11] the analogs of
Morse inequalities were extracted. These inequalities
quickly became an object of study and applications.
Farber [3] obtained an exactness theorem
for these inequalities in the case
$\pi_1 (M)=\ZZZ~,~ \dim M\geq 6$. The author
 [14, 15] obtained an exactness result for these inequalities
in the case $\pi_1 M=\ZZZ^m,\dim M\geq 6$.
These inequalities were applied by J.Cl.Sikorav [22]
to the theory of Lagrangian intersections.

The study of the properties of the chain complex itself
advanced more slowly. In [21] a refined version
of the complex, defined over a completion of the group
ring of $\pi_1(M)$ was defined. In [16, 17] I proved that
this complex computes the completed simple homotopy type
of $M$ and in [18, 19] I obtained an exactness result
for the case of general fundamental group. This theorem
implies the results of [3] and [14]. The algebraic
result of Ranicki [20] show that this result implies also
the Farrell fibration theorem [4]. (Note that the
resulting proof of the Farrell's result is independant
of Farrell's original proof).

 One of the objects here, which behaves
very differently from the classical situation
is
the boundary
operator in the Novikov complex.
Fix some $k$. The boundary operator
 $\partial :C_k(f)\to C_{k-1}(f)$
is represented by a matrix, which entries are in
the ring of Laurent power series. That is
$\partial_{ij}=\sum_{n=-N}^{\infty} a_nt^n,\wwhere
 a_n\in\ZZZ$.

Since the beginning S.P.Novikov conjectured that the
power series
$\partial_{ij}$
 have some nice analytic
properties. In particular he conjectured that
\vskip0.1in
\it
Generically the coefficients $a_n$ of
$\partial_{ij}=\sum_{n=-N}^{\infty} a_nt^n$
grow at most exponentially with $n$.
\rm
\vskip0.1in
See [13, p. 229] and [1, p.83]
for (different) published versions
of this conjecture.
In the present paper we prove that generically
the incidence coefficients are rational functions
of $t$ of the form
$\frac{P(t)}{t^k Q(t)} \wwhere
P(t),Q(t)$
are polynomials, $k\in \NNN$ and $Q(0)=1$.

The methods developed in the present paper
allow to make explicit computations of the Novikov
complex. These computations and
further generalizations are the subject
of the second part of the work, to appear.

We would like to conclude this subsection by
the following remark.
      Studying the algebraic
properties of the Novikov complex
and of the Novikov completion
one passes often from
localization to completion
and vice versa (see [3], [13, 14],
[20]). We see now that these passages have a geometric
background.
\subhead{B. Statement of the main theorem}\endsubhead
Let $M$ be a closed connected manifold and $f:M\to S^1$
a Morse map, non homotopic to zero. Denote the set of
critical points of $f$ by $S(f)$.
The set of  $f$-gradients
of the class  $\smo$,
 satisfying the transversality
assumption (see \S1 for terminology),
 will be denoted by ${\Cal G} t(f)$.
By Kupka-Smale theorem it is residual in the set of
all the $\smo$ gradients. Choose
$v\in\gt (f)$. Denote by $\bar M @> {\Cal P} >>  M$
the connected infinite cyclic covering for which
$f\circ \Cal P$ is homotopic to zero.
Choose a lifting  $F:\bar M\to\RRR$ of
$f\circ \Cal P$ and let $t$ be the generator
of the structure group of $\Cal P$
such that $F(xt)<F (x)$.
The $t$-invariant lifting of $v$ to $\bar M$
will be denoted by the same letter $v$.
For every critical point $x$ of $f$ choose a
 lifiting $\bar x$
of $x$ to $\bar M$. Choose orientations
of stable manifolds of critical points.
Then for every $x,y\in S(f)$ with $\ind x =\ind y +1$
and every $k\in\ZZZ$ the incidence coefficient
$n_k(x,y;v)$
is defined
(as the algebraic number of $(-v)$-trajectories
joining $\bar x$ to $\bar y t^k$).

\proclaim{Main Theorem}
In the set $\gt (f)$ there is a subset $\gt _0(f)$ with
the following properties:
\roster\item
$\gt_0 (f)$ is open and dense in $\gt (f)$ with respect to
$C^0$ topology.
\item If $v\in\gt _0(f),~x,y\in S(f)\aand \ind x=\ind y+1$,
then
$\sum_{k\in\ZZZ} n_k(x,y;v)t^k$ is a rational
function of $t$ of the form
$\frac{P(t)}{t^m Q(t)}, \wwhere
 P(t) \aand Q(t)$
are polynomials with integral
coefficients,
$m\in\NNN$,~and~$Q(0)=1$.
\item Let $v\in\gt _0(f)$.
Let $U$ be a neighborhood of $S(f)$.
Then for every $w\in\gt _0(f)$ such that
$w=v$ in $U$ and $w$
is sufficiently close to $v$ in $C^0$ topology
we have:
$n_k(x,y;v)=n_k(x,y;w)$
for every $x,y\in S(f),~k\in\ZZZ$.
\endroster
\endproclaim

\remark{Remark} This theorem implies easily
the exponential
estimate from the Novikov conjecture above:
note that every rational function of the form
$\frac {P(t)}{Q(t)}$, where $Q(0)\not= 0$,
has a non zero radius of convergency in $0$.
  \endremark
\subhead C. Further generalizations \endsubhead
In the second part of this work, to appear, we
develop our methods and generalize them to
the incidence coefficients with values
in the Novikov completions of the group rings.
We prove that for a $C^0$ generic $f$-gradient these
coefficients belong actually to the image in the
Novikov ring of the corresponding localization of the group ring
of the fundamental group. For the case of
irrational Morse forms we prove the analog of this result,
concerning the incidence coefficients related with free abelian
coverings. We state and prove the generalization of
the exponential growth rate estimate for the case
of  Morse forms of arbitrary
irrationality degree (and the incidence coefficients related to the universal
covering). We also give an example of explicit computation
of the Novikov complex.

\subhead D. Remarks on the contents of the paper \endsubhead
The proof of the main theorem is contained in
\S 4.E. The technical results on Morse
functions and their gradients
which we use are that of [19,\S 2].
We reformulate them (in a slightly
generalized and modified form)
in \S 1.
\S 5 contains some results on $C^0$ stability
of trajectories of $C^1$ vector fields. These results
can be of independent interest.
\S 5 is independent of other sections
and the results of \S 5 are used in \S 1 - 4.

\subhead E. Remarks on the terminology \endsubhead
If $W$ is a manifold with boundary, then
$\Wkr$ stands for $W\setminus \partial W$.
$\gamma(x,t;v)$ stands always for
the value at $t$ of the integral curve
of the vector field $v$, starting at $x$.
If $\gamma(x,\cdot ;v)$ is defined
on $[\alpha, \beta]$, then
$\{\gamma(x,t;v)\mid t\in[\alpha,\beta]\}$
is denoted by $\gamma(x,[\alpha,\beta];v)$.
For a (time dependant) vector field
$w$ on a closed manifold $M$
and $t\in\RRR$ we denote by $\Phi(w,t)$
the diffeomorphism $x\mapsto\gamma(x,t;w)$
of $M$.
The overline denotes closure (e.g. $\overline U$).
The bar is reserved for the objects related to the
infinite cyclic covering (e.g. $\Cal P :\bar M\to M$).

The given riemannian metric on
a riemannian manifold $M$ will be
denoted by $\rho$, the corresponding
 norm on $T_xM$ will be
denoted by
$\mid\cdot\mid_\rho$, the $C^0$ norm
 on the space of vector fields on $M$
 will be denoted by
$\mid\mid\cdot\mid\mid_\rho$.
The euclidean  metric in ${\bold R}^n$
 will be denoted by
$e$. When there is
no possibility of confusion we drop
the indices $\rho$ and $e$
 from the notation.

For a chart $\Phi:U\to V\subset\RRR^n$
of $M$ and $x\in U$
we denote by $\GG(x,\Phi)$
the number
$\sup\limits_{h\in T_xM, h\not=0}
(\max \big(\vert h\vert_\rho/\vert\Phi_*h\vert_e,
\vert\Phi_*h\vert_e/\vert h\vert_\rho)\big)$.
 For $K\subset U$
we denote by $\GG(K,\Phi)$ the number
$\sup_x
        \GG(x,\Phi)$, where $x$ ranges over points of $K$.

The end of a definition, remark or construction
is marked by $\triangle$.

\subhead  F. Acknowledgements \endsubhead
I am grateful to S.P.Novikov, J.Prszytycki, J.-Cl.Sikorav,
\break
C.Simpson, T. tom Dieck,
A.N.Tyurin, P.Vogel
for valuable discussions.

A part of this work was done during my stay
in G\"ottingen University
in the spring 1994. It is a pleasure
to express here my gratitude
 to Sonderforschungsbereich
170 Analysis and Geometry
for hospitality and financial support.

\head \S1 Preliminaries on Morse functions
and their gradients \endhead
Subsections A-D contain generalities. In
Subsection G we formulate Theorem 1.17,
which is the main aim of \S\S
 1 and 3.
 Subsections E,F
contain some techniques, useful for its proof.
All the objects
are supposed to be of class $\smo$.

\subhead A. Terminology: functions and gradients
 \endsubhead
In the present and the following
subsections $M$ stands for a closed manifold.
\definition{Definition 1.1} Let $f:M\to \RRR$ be
 a Morse function on a closed manifold
$M$. Denote $\dim M $ by $ n$.
The set of critical points of $f$ will be
denoted by $S(f)$.
A chart $\Phi_p:U_p\to B^n(0,r_p)$
(where $p\in S(f), U_p$ is a neighborhood
of $p$, $r_p >0$) is called \it standard
chart for $f$ around $p$ of radius $r_p$
\rm (or simply $f$-\it chart\rm ) if
there is an extension of $\Phi_p$ to a chart
$\widetilde\Phi_p:V_p\to B^n(0,r'_p)$,~
(where $\overline {U_p}\subset V_p\aand r'_p>r_p$),
such that
$(f\circ \widetilde\Phi_p^{-1})(x_1,...,x_n)=f(p)+
\sum_{i=1}^n\alpha_i x_i^2$, where
 $\alpha_i<0$ for $i\leq \ind_fp$ and
$\alpha_i>0$ for $i>\ind_fp$.
 The domain $U_p$
is called \it standard coordinate neighborhood.
\rm
Any such extension $\widetilde\Phi_p$ of $\Phi_p$
will be called
\it standard extension of $\Phi_p$.
\rm
The set $\Phi_p^{-1}(\RRR^k\times\{0\})$, resp.
$\Phi_p^{-1}(\{0\}\times\RRR^{n-k})$,
where $k$ stands for $\ind_fp$,
is called \it negative disc, \rm
resp. \it positive disc. \rm
If for every $i$ we have $\alpha_i=\pm 1$,
we shall say that the coordinate system
$\{\Phi_p\}$ is
\it strongly standard. \rm

A family
${\Cal U}= \{\Phi_p:U_p\to B^n(0,r_p)\}_{p\in S(f)}$
 of $f$-charts is called $f$-\it chart-system\rm ,
 if the family $\{\overline{U_p}\}$ is disjoint.
We denote $\min_p r_p$ by $d({\Cal U})$, and
 $\max_p r_p$ by $D({\Cal U})$.
If all the $r_p$ are equal to $r$, we shall say
that $\UU$
\it is of radius $r$. \rm
The set $\Phi_p^{-1}(B^n(0,\lambda))$, where
 $\lambda\leq r_p$ will be denoted by $U_p(\lambda)$.
For $\lambda\leq d(\UU)$ we denote
$\cup_{p\in S(f)}U_p(\lambda)$ by $\UU(\lambda)$.

Let $\UU=\{\Phi_p:U_p\to B^n(0,r_p)\}_{p\in S(f)}$,
$\UU'=\{\Phi'_p:U'_p\to B^n(0,r'_p)\}_{p\in S(f)}$
be two $f$-chart-systems. We say, that $\UU'$
is
\it
a restriction of
\rm
 $\UU$, if for every $p\in S(f)$
we have:
 $r'_p\leq r_p, ~U'_p\subset U_p, \Phi'_p=\Phi_p\mid U'_p$.

Given an $f$-chart system
${\Cal U}= \{\Phi_p:U_p\to B^n(0,r_p)\}_{p\in S(f)}$
we say, that a vector field $v$ on $M$
 is an
\it $f$-gradient with respect to $\Cal U$,
\rm
if

1) $\forall x\in M \setminus S(f)$ we
have $df(v)(x)>0$;

2) $\forall p\in S(f)$ we have
$(\widetilde\Phi_p)_*(v)=(-x_1,...,-x_k,x_{k+1},...,x_n)$,
where $k=\ind_f p$, and $\widetilde\Phi_p$ is some standard
extension of $\Phi_p$.

We say that a vector field $v$ is an
\it $f$-gradient \rm
if there is an $f$-chart system
 $\Cal U$ ,
such that $v$ is an $f$-gradient
with respect to $\Cal U$.
\quad
$\triangle$
\enddefinition

\definition{Definition 1.2}
Assume that $M$ is riemannian.
Let $\delta_0>0$ be the radius of injectivity of $M$.
For $0<r<\delta_0$ and $x\in M$
we denote by $D_r(x)$, resp. $B_r(x)$,
the image of  $D^n(0,r)$, resp. $B^n(0,r)$,
with respect to the exponential map $\text{exp}_x$.
Let
$f:M\to\RRR$ be a Morse function,
$v$ be an $f$-gradient, $p\in S(f)$. Set:
$$\gather
B_\delta(p,v)=\{x\in M~\vert~ \exists
t\geq 0 :
 \gamma (x,t;v)\subset B_\delta (p)\}\\
D_\delta(p,v)=\{x\in M~\vert ~\exists
t\geq 0 :
\gamma (x,t;v)\subset D_\delta (p)\}~~,~~
D(p,v)=\{x\in M\vert \lim\limits_{t\to\infty}
\gamma(x,t;v)=p\}
\endgather
$$
For $s\in \NNN, 0\leq s\leq n$
we denote by $B_\delta (\indl s;v)$
the union of   $
B_\delta (p,v)$ where $p$ ranges over critical
points of $f$ of index $\leq s$.
 Similar notations like
$D_\delta (\indl s;v)$
or $K(\inde s;v)$ etc. are now
clear without special definition.
\quad $\triangle$\enddefinition

\remark{Remarks}
1) Our definition of $f$-gradient is
wider than that of [18, 19] (and than that of
gradient-like vector field in [6]),
since we allow $\alpha_i\not=\pm 1$
in 1.1. It is not difficult to show
that all
 the results from [19,\S 2]  rest true
with the present definition.
The terminology of [19,\S2] concerning
Morse functions and their gradients
is carried over to the present case
without changes.
For convenience of the reader we include the results
and some of the terminology from [19,\S2]
which we use here,
to Subsection C.

2) One can define these sets for any $\delta>0$, see [19,\S2];
they were denoted there by
$B_\delta(\leq s;v)$, resp.
$D_\delta(\leq s;v)$.

3) Similarly to 1.1 one defines the notion of
$f$-gradient for Morse maps $f:M\to S^1.\triangle$

\endremark

\subhead
B. Terminology: $\MM$- flows
\endsubhead
\definition{Definition 1.3}
Assume that $M$ is riemannian.
 We say, that a triple $\flow$ is
an $\MM$-flow on $M$
($\MM$ for Morse) if
(1) $f:M\to\bold R$ is a Morse function,\quad
(2) $\UU$ is an $f$-chart system, \quad
(3) $v$ is an $f$-gradient with respect to $\UU$,\quad
(4) for each $f$-chart $\chart$ the coordinate
 frame in $p$ is orthonormal with respect to
 the riemannian metric of $M$.\quad$\triangle$
\enddefinition

Let $\flow$ be an $\MM$-flow on $M$,
where
$\UU=           \atlas$.

We       denote
       $\max_p\GG(U_p,\Phi_p)$ by
$\GG(\VV)$.
 We       say that $\VV$
\it is
of radius $r$, \rm
if~  $\UU$ is of radius $r$.

\vskip0.1in
 Let $\gamma$ be a $v$-trajectory.

A) The number
of sets $U_p=U_p(r_p)$,
intersected by $\gamma$, will be denoted by
$N(\gamma)$. The number $\max_\gamma N(\gamma)$
will be denoted by $N(\VV)$.
The set
$\{t\in\RRR \vert \gamma(t)\notin\cup_{p\in S(f)}
U_p\}$ is a finite union of closed intervals; its measure
 will be denoted by $T(\gamma)$.
Let $\beta>0,~C>0$. We say, that $\VV$ is
\it
$(C,\beta)$-quick,
\rm
if
$\nrv \leq C$
and for every $v$-trajectory $\gamma$
we have $T(\gamma)\leq\beta$.

B) Let
$\delta\leq d(\UU)$.
The number
of sets $U_p(\delta)$,
intersected by $\gamma$, will be denoted by
$N(\gamma,\delta)$. The number
$\max_\gamma N(\gamma,\delta)$
will be denoted by $N(\VV,\delta)$.
The set
$\{t\in\RRR \vert \gamma(t)\notin\UU(\delta)\}$
 is a finite union of closed intervals
and its measure
 will be denoted by $T(\gamma,\delta)$.
Let $\beta>0,~C>0$. We say, that $\VV$ is
\it
$(C,\beta,\delta)$-quick,
\rm
if
$\nrv \leq C$
and for every $v$-trajectory $\gamma$
we have $T(\gamma,\delta)\leq\beta$.

\vskip0.1in
Let  $\VV^{(i)}=(f^{(i)},v^{(i)},\UU^{(i)})$, where i=1,2 be two
$\MM$-flows on $M$. Set
$\UU^{(i)}= \{\Phi_p^{(i)}:U_p^{(i)}\to B^n(0,r^{(i)}_p)\}_{p\in S(f^{(i)})}$.
 We say, that $\VV^{(2)}$
\it is subordinate to \rm $\VV^{(1)}$, if
(1) $\UU ^{(2)}$ is a restriction of $\UU^{(1)}$,\quad
(2) $v^{(2)}=\phi\cdot v^{(1)}$, where
$\phi :M\to\bold R^+$ is a $\smo$
 function such that $\phi(x)=1$ for $x$
 in a neighborhood of the closure of
$\cup_pU_p^{(2)}$.

\remark{Remark} There is an obvious analog of the
  terminology of Subsections A, B
           for a more general situation of
Morse functions on
compact
 cobordisms. We shall make free use of it.
~$\triangle$
\endremark

\subhead{C. First properties of $f$-gradients}\endsubhead
In this subsection $f:W\to [a,b]$ is a Morse function on a compact
cobordism,  $V_0=f^{-1}(a),~V_1=f^{-1}(b)$,~$v$ is an $f$-gradient.
Denote $\cup_{p\in S(f)}B_\delta(p,v)$ by $B_\delta(v)$,
$\cup_{p\in S(f)}D_\delta(p,v)$ by $D_\delta(v)$,
$\cup_{p\in S(f)}D(p,v)$ by $K(v)$.
Denote $\dim W$ by $n$.
\proclaim{Lemma 1.4}
\roster\runinitem For every $p\in S(f)$ the set $B_\delta(p,v)$ is open.
\item $D_\delta(v)$ and $K(v)$ are compact.
\item $D_\delta(v)=\cap_{\theta>\delta} B_\theta (v)\aand
K(v)=\cap_{\theta>0}B_\theta (v)$.
\item $D_\delta(v)=\overline{B_\delta(v)}$.
\endroster
\endproclaim
\demo{Proof} Proof of (1) - (3) is similar to that of
Lemma 2.3 of [19].
(4) follows from (2) and the (obvious) inclusions
$B_\delta(v)\subset D_\delta(v)\subset\overline{B_\delta(v)}$.
\quad$\square$            \enddemo

We shall accept here the terminology of [19,\S2].
We only recall from there that
 $v$ is called:
\roster\runinitem
      {\it perfect }, if $\big( x,y \in S(f) \big)
\Rightarrow \big( D(x,v) \capp D(y,-v) \big)$
      \item
      {\it almost good}, if $\big( x,y \in S(f)$ and
$\ind x \leq \ind y \big) \Rightarrow \big(
D(x,v) \capp D(y,-v)\big)$.
    \endroster
\vskip0.1in
We say that a $f$-gradient is
\it
$\delta$-separated
\rm
if there is an ordered Morse function
$\phi:W\to[a,b]$, adjusted to $(f,v)$,
with an ordering sequence
$a_0,...,a_{n+1}$, such that
\roster\item for
every $p\in S(f)$ we have $D_\delta(p)\subset
\phi^{-1}(]a_k,a_{k+1}[)$, where $k=\ind p$.
\item for every $l:0\leq l\leq n$ and every
$q\in S(f), \ind q=l$
there is a Morse function
$\psi:W_l\to[a_l,a_{l+1}]$, where
$W_l=\phi^{-1}([a_l,a_{l+1}])$,
adjusted to $(\phi\vert W_l,v\vert W_l)$
and a regular value $\mu\in]a_l,a_{l+1}[$ of $\psi$, such that
$D_\delta(q)\subset\phi^{-1}(]a_l,\mu[)$
and for every $r\in S(f),\ind r= l, r\not= q$ we have
$D_\delta(r)\subset\phi^{-1}(]\mu,a_{l+1}[)$.
                         \endroster

Note that if $v$ is $\delta$-separated, we have
$D_\delta(p,v)\cap D_\delta(q,-v)=\emptyset$
 whenever
$\ind p\leq \ind q$.
The proofs of the next two results
are similar to that of
[19, Le.   2.7] and [19,Le.   2.11].
\proclaim{Lemma 1.5} If $v$ is almost good,
then $v$
is $\delta$-separated for
some $\delta>0$.
~~$\square$\endproclaim

\proclaim{Proposition 1.6}
If $v$ is  $\delta_0$-separated, then
$\forall
\delta\in]0,\delta_0[$ and
$\forall
s: 0\leq s\leq n$
$$\aligned
&(1)~~D_\delta (\indl s;v)~ \text{~~is compact}~~;\\
&(3)~~{\bigcap} _     {\theta>0} B_\theta (\indl s;v) =
K (\indl s;v);\\
&(5)~~\overline{B_\delta (\indl s;v)} = D_\delta (\indl s;v).
\endaligned\qquad
\aligned
&(2)~~K (\indl s;v) ~ \text{~~is compact}~~;\\
&(4)~~{\bigcap} _{\theta>\delta} B_\theta (\indl s;v)
= D_\delta (\indl s;v);\\
&\qquad\quad\qquad\square
\endaligned
$$
\endproclaim

For two regular values
$\lambda,\mu$ of $f$
with $\lambda<\mu$ we denote by
$K^+_\mu$ the set $f^{-1}(\mu)\cap
\bigcup\limits_{p\in S(f)\cap f^{-1}([\lambda,\mu])}
D(p,-v)$;
 by $K_\lambda^-$
the set
$f^{-1}(\lambda)\cap \bigcup\limits_{p\in S(f)
\cap f^{-1}([\lambda,\mu])}
D(p,v)$
and by $\stind v\mu\lambda$ the
$\smo$ diffeomorphism
$f^{-1}(\mu)\setminus K_\mu^+\to f^{-1}(\lambda)
\setminus K_\lambda^-$,
which associates to each point $x$ the point of
 intersection
of $\gamma(x,\cdot;-v)$ with $f^{-1}(\lambda)$.
For $X\subset f^{-1}(\mu)$
we denote (by abuse of notation)
$\stind v\mu\lambda (X\setminus K_\mu^+)$
by $\stind v\mu\lambda (X)$.
If the values of $\mu$ and $\lambda$
are clear from
the context we suppress them in the notation.

\subhead{D. $C^0$-stability properties }
\endsubhead
In this subsection $f:W\to[a,b]$ is a Morse
function on a compact riemannian
 cobordism $W$,
$f^{-1}(a)=V_0, f^{-1}(b)=V_1$,
$v$ is an $f$-gradient.

\proclaim{Proposition 1.7}
Let $\delta>0$. Let $K$ be a compact
in $V_0\cap B_\delta(v), R_1$ be an open neighborhood
of $D_\delta(v)\cap V_0, R_2$ be an open
neighborhood of $K(v)\cap V_0$.
Then there is $\epsilon>0$ such that
for every $f$-gradient
$w$ with $\Vert w-v\Vert<\epsilon$
we have:

\text{\rm (1)}~~~
       $K\subset B_\delta(w)$,\quad
\text{\rm (2)}~~~
     $D_\delta(w)\cap V_0\subset R_1$,\quad
\text{\rm (3)}~~~
      $K(w)\cap V_0\subset R_2$.
\endproclaim
\demo{Proof}
(1) is proved by a compactness argument
similar to the proof of 5.6.
To prove (2) note that
 $V_0\setminus R_1$ is a compact such that
each $(-v)$-trajectory
starting at a point of  $V_0\setminus R_1$
reaches the boundary and that
$\tau(V_0\setminus R_1,-v)\subset W\setminus
\cup_{p\in S(f)} D_\delta (p)$.
By 5.6 the same is true for every $w$ sufficiently close to
$v$ in $C^0$ topology.
The proof of (3) is similar.
~$\square$
\enddemo

Next lemma is obvious.
\proclaim{Lemma 1.8}
Let $g: W\to[a,b]$ be a Morse function, adjusted to $(f,v)$.
Then there is $\epsilon>0$ such that
every $f$-gradient $w$ with
$\Vert w-v\Vert <\epsilon$
is a $g$-gradient.\quad$\square$
\endproclaim

\proclaim{Corollary 1.9}
If $v$ is almost good, resp.
$\delta$-separated,
then there is $\epsilon>0$
such that every $f$-gradient $w$ with
$\Vert w-v\Vert<\epsilon$
is almost good,
resp.
$\delta$-separated.
\quad$\square$\endproclaim

\proclaim{Proposition 1.10}
Let $p\in S(f)$. Let $U$
be an open set in $V_0$. Assume that for every $x\in D(p,v)$
the trajectory $\gamma(x,\cdot;-v)$
reaches the boundary and that
$D(p,v)\cap V_0\subset U$.
Then there is $\epsilon>0$ such that for every
$f$-gradient $w$ with $\Vert w-v\Vert<\epsilon$
and every $x\in D(p,w)$
the trajectory $\gamma(x,\cdot;-w)$ reaches the
boundary and $D(p,w)\cap V_0\subset U$.
\endproclaim
\demo{Proof}
By Rearrangement theorem
[19, Th. 2.6]
 there is a Morse function
$\phi:W\to[a,b]$ adjusted to $(f,v)$, and
a regular value $\mu$ of $\phi$ such that
$\phi(p)<\mu<\phi(q)$
for every $q\in S(f), q\not= p$.
Now apply to $\phi\vert\phi^{-1}([a,\mu])$ the previous
lemma and 1.7(3).
\quad$\square$
\enddemo

\subhead E. Two lemmas on standard gradients
 in $\bold R^n$ \endsubhead
During this subsection we refer to $\bold R^n$ as to
the product $\RRR ^k\times \RRR ^{n-k}$;
a point $z\in \RRR ^k$ is therefore denoted
by $(x,y)$, where $x\in \bold R^k,y\in \bold R^{n-k}$.
The vector field with components $(x,-y)$ is denoted
by $v_0$; the standard euclidean norm in $\RRR ^n$
will be denoted by $\mid \cdot \mid$.  The Morse
function $(x,y) \mapsto -\mid x\mid^2 + \mid y \mid ^2$
 will be denoted by $f_0$.
(Then
       $-2v_0$ is the riemannian
gradient of $f_0$ with respect to the euclidean metric.)
The $v_0$-trajectories are all of the form
$(x_0e^t,y_0e^{-t})$. Using this fact,
 the following two lemmas become
the matter of computation, which
will be left to the reader.

\proclaim{Lemma 1.11}
Let $R>r>0$ and $\gamma$ be a
$v_0$-trajectory. Then the time,
which $\gamma$ spends in $B^n(0,R) \setminus
B^n(0,r)$ is not more than
$\ln \big( (\frac Rr)^2
 +\sqrt{(\frac Rr)^4 -1}~\big) $
 and the length
 of the corresponding part of $\gamma$ is
 not more than $2R$.\quad$\square$
\endproclaim
\proclaim{Lemma 1.12}
Let $r>0$ and $\gamma$ be a $v_0$-trajectory.
 Then the time, which $\gamma$ spends in the set
$f_0^{-1}([-r^2,r^2]) \setminus B^n(0,r)$ is
not more than $2$.\quad$\square$
\endproclaim

\proclaim{Corollary 1.13}
Let $\flow$ be an $\MM$-flow on
a closed manifold $M$. Assume, that
$\VV$ is $(C,\beta,\alpha)$-quick.
Then
 $\VV$ is $(C,\beta + 8N(\VV),\alpha /2)$-quick.
\endproclaim
\demo{Proof}
By 1.11
the time, which a $v$-trajectory can
 spend in
$\cup_{p\in S(f)} (U_p(\alpha)\setminus U_p(\alpha /2))$
is not more than
$N(\VV )\ln(2^2+\sqrt{2^4 -1})\leq 8N(\VV )$.
$\square$
\enddemo
\remark{Remark 1.14} Corollary 1.13
is true as it stands for Morse functions on
cobordisms.\quad $\triangle$
\endremark

\subhead F. A-construction \endsubhead
Let $\flow$ be an $\MM$-flow on a closed
manifold $M$,
where $\UU =
\{\Phi_p:U_p\to B^n(0,r_p)\}_
{\scriptscriptstyle {p\in S(f)}}$
and $n=\dim M$.
 Let $\mu<\min_pr_p$.
 Denote $U_p(\mu)$ by $U'_p$.
 Let $\Gamma \geq1$.
We shall construct an
$\MM$-flow $\VV ' = (f,w,\UU ')$, subordinate
 to $\VV$ where $w=\phi\cdot v$, and
$\phi:M\to \RRR ^+, \phi (x) =1$ for
$x\in \obb '$ and $\phi (x) =\Gamma$
for $x\notin \obb $. The main property
of the gradient $w$ is an explicit estimate of
 the time, which a $w$-trajectory can
 spend inside each $U_p\setminus U'_p$.
 This construction will be used in the subsection G.
Since we have set $w(x)=\Gamma\cdot v(x)$
for $x\notin \obb $ we have only to define
 $w(x)$ inside each $U_p$ for $\crr$, which
will be done in 1.15.
Denote $\nrv$ by $B$, and
$\Cal G (\VV)$ by $D$.

\definition{ Construction 1.15}
In 1.15 we deal only with
a neighborhood of one critical
 point (say, $p$),
so we drop the index $p$ from
the notation. The index
of this critical point is denoted by $k$,
the standard chart
is denoted by
$\Phi:U\to B^n(0,r)$ and we are
given $\mu$ with $0<\mu <r$.
Denote the function $x\mapsto
\ln (x^2 +\sqrt{x^4 - 1})$ by $\LLN (x)$.
Choose $\delta\in ]0,\frac{r -\mu}{2} [$
 so small that
$$
\LLN\big(r/(r-\delta)\big)<\min(\frac {Dr}{2B} ,1/2),
\quad
\LLN\big((\mu+\delta)/\mu\big)<\min(\frac {Dr}{2B} ,1/2)
\tag1.1
$$

Let $\theta :\bold R\to [0,1]$ be a $\smo$
function such that
$\supp\theta\subset ]\mu,r[$ and
 $\theta (x)=1$ for $x\in [\mu+\delta,r-\delta]$.
Let $\lambda :[0,r]\to\RRR^+$ be the
$\smo$ function, defined by the following
formulas:
$$\gather
\lambda (t)=(1-\theta (t))+\theta(t)\frac B{Dt}
 \ffor t\in [0,\mu+\delta];\qquad \lambda(t)=
\frac B{Dt}\ffor\in [\mu+\delta,r-\delta]\\
 \lambda(t)=\Gamma(1-\theta(t))+\theta(t)\frac B{Dt} \ffor
t\in [r-\delta,r]
\endgather
$$

Note that $\lambda (t) =1$ for $t\in [0,\mu],~\lambda (t) =
\Gamma$ nearby $r$.
Define a vector field $v_1$ in $B^n(0,r)$ by
$v_1(x)=\lambda(\vert x\vert_e)v_0(x)$, where
 $v_0(x)=(-x_1,...,-x_k,x_{k+1},...,x_n)$.
Note that $(\Phi_*^{-1})(v_1)$
 equals $\Gamma v$ in
$\Phi^{-1}(B^n(0,r)\setminus B^n(0,r-\nu))$
for some small $\nu>0$.
\quad$\triangle$
\enddefinition

Applied to the neighborhood $U_p$
of each critical point $p$ this construction
 extends $w$  to the whole of $M$.
Denote by $\Phi_p'$ the restriction of
 $\Phi_p$ to $U'_p$ and set
$\UU'=\{\Phi_p'\}_{\crr}$.
It is obvious from the construction that
 $\VV ' = (f,w,\UU ')$ is an $\MM$-flow
on $M$,
subordinate to $\VV$.
 We shall denote it by
$(A)(\VV , \Gamma,\mu)$.

\proclaim{Lemma 1.16}
\roster
\runinitem $w(x) =\Gamma v(x)$ for
$x\notin \obb \aand\nrw _\rho\leq \Gamma\nrv _\rho$;
\item  For any $w$-trajectory $\gamma$ and
every $\crr$ the time, which $\gamma$
 spends inside
$U_p\setminus U'_p$ is not more, than
$\big(3\GG (\VV)r_p\big)/\Vert v\Vert_\rho$.
\endroster
\endproclaim

\demo{Proof} The first part of (1) is obvious
from the construction.
The estimate of $\Vert w\Vert_\rho$
is obtained by an explicit computation
in the standard coordinate systems
which we leave to the reader.
 To prove (2) note that the
time which
a $v_1$-trajectory can spend in
$B^n(0,r)\setminus B^n(0,r-\delta)$
or in
$ B^n(0,\mu+\delta)\setminus B^n(0,\mu)$
is not more than $\frac {Dr}{2B}$, which follows directly
from the definition of $v_1$, Lemma 1.11 and
(1.1).
Further, the euclidean
length of the part of the $(v_1)$-trajectory
in $\big( B^n(0,r -\delta)\setminus
B^n(0,\mu +\delta)\big)$ is not more
than $2r$ (by Lemma 1.11), and the euclidean
norm of the tangent vector is equal to $B/D$.
Therefore the time spent in
$B^n(0,r-\delta)\setminus
B^n(0,\mu+\delta)$ is not more than $(2Dr)/B$,
and the total time spent in
$B^n(0,r)$ is not more than $(3Dr)/B$.\quad$\square$
\enddemo

\subhead{G. Statement of the theorem on the existence
of quick flows}\endsubhead
\proclaim{Theorem 1.17} Let $M$ be a
closed riemannian manifold, $\dim M = n$.
Let $C>0,
\mathbreak
\beta>0, A>1$.
Then there is an almost
 good $\MM$-flow $\flow$ of radius $r$,
  and for every
$\mu\in ]0, r[$
there is an
$\MM$-flow $\WW = (f,w,\UU ')$,
subordinate to $\VV$, and such that

\text{\rm (1)}\quad $\WW$ is of radius $\mu$ and
 is $(C,\beta )$-quick;\qquad
\text{\rm (2)}\quad $N(\WW  ) \leq n2^n$;\qquad
\text{\rm (3)}\quad $\GG (\WW )\leq A$.

\endproclaim

This theorem will be proved by induction in
$n$ with the help of the theorem 1.18.

\proclaim{Theorem 1.18} Let $M$ be a closed
 riemannian manifold, $\dim M =n$. Let $B>0,\mu_0>0$.
Then there is an almost good $\MM$-flow $\flow$
on $M$ with the following properties:

\text{\rm (1)}\quad $\GG (\VV )\leq 2\aand D(\VV ) <\mu_0 $;\qquad\qquad
\text{\rm (2)}\quad$N(\VV )\leq n2^n$;

\text{\rm (3)}\quad$\VV$ is $(B,R(n))$-quick,
where $R(n) = 100+n2^{n+7}$.
\endproclaim

The following corollary of 1.17 will be used in \S 4.

\proclaim{Corollary 1.19}
Let $M$ be a closed riemannian manifold,
$\dim M =n$. Let $C>0,\beta>0$.
Then there is an almost good
$\MM$-flow $\flow$ on $M$, such that
 for every $\delta>0$
sufficiently small there is an $\MM$-flow
$\WW =(f,w,\UU ')$, subordinate to $\VV$
and such that \roster
\item
$\nrw\leq C$ and $\UU'$ is of radius $\delta/2$;
\item For every $t\geq \beta$ and every
$s:0\leq s\leq n$ we have
$$\gather
\Phi (-w,t)\big(M\setminus B_\delta (\indl {n-s-1};-v)\big)
 \subset B_\delta (\indl s;v)\tag 1.2\\
\Phi (w,t)\big(M\setminus B_\delta (\indl {n-s-1};v)\big)
 \subset B_\delta (\indl s;-v)\tag 1.3           \endgather
$$
\endroster
\endproclaim
\demo{Proof}
Let $\VV = (f,v,\UU)$ be the
$\MM$-flow, satisfying 1.17 with respect
 to $C,\beta$, and $A=2$. Let $\delta>0$
be so small that
$\forall \crr$ the disc $\rho(p,x) \leq\delta$
belongs to $U_p$. Since $A=2$,
 that implies, that this
disc contains $U_p(\delta/2)$.
 Let $\WW = (f,w,\UU ')$ be the
  $\MM$-flow,
satisfying (1)-(3) of 1.17
 with respect to $\delta/2$. I claim,
that $\WW$ satisfies our conclusion.
 Indeed, 1) goes by construction.
To prove 2), let
$y\in M\setminus B_\delta (\indl {n-s-1};-v)$.
Since $\WW$ is $(C,\beta)$-quick,
there is $t_0\in [0,\beta]$, such that
$\gamma (y,t_0;-w)\in U_p(\delta/2)$ for some
 $p\in S(f)$.
Since $M\setminus B_\delta (\indl {n-s-1} ;-v)$
is $v$- and $w$-invariant, we have
 $\ind p\leq s$.
Since $B_\delta (\indl s; v)$ is $(-v)$- and $(-w)$-
invariant, we have
$\gamma (y,\lambda;-w)\in B_\delta (\indl s;v)$
for all $\lambda\geq t_0$.
Proof of (1.3) is similar.
\quad$\square$

\enddemo
The scheme of the proof of 1.17 and 1.18 is
as follows. The proof of 1.18 $\Rightarrow$
1.17 is given below. The proof of
 1.17(n) $\Rightarrow$  1.18(n+1)
 is done in \S3 with the help of what we call
(S)-construction. This proof also gives the proof
of 1.18(1).

\demo{ Proof of 1.18 $\Rightarrow$ 1.17}
We are given $C,\beta$ and $A$. Choose $C'$
such that $C'\leq C$
and that $R(n)\frac{C'}C\leq\beta/2$.
Then choose $\mu_0>0$ such that
 $6n2^n \frac{\mu_0}{C'} < \beta/2$.
Let $\flow$ be an $\MM$-flow on $M$, which
satisfies the conclusions of 1.18 with respect
to $\mu_0$ and $C'$.
 Here $\UU=\atlas$.
Choose $r>0$ so small that $r<\min_p r_p$
and $\GG (p,r)\leq A$ for every $p\in S(f)$.
Denote by $\UU'=\{\Phi'_p:U_p'\to
B^n(0,r)\}_{p\in S(f)}$
the corresponding restriction of $\UU$.
We claim that the $\MM$-flow
$\VV '=(f,v,\UU')$ satisfies the conclusions of 1.17.
Indeed, let $0<\mu<r$.
 Apply the $A$-construction to
$\VV',~\Gamma=C/C'$ and to
$\mu$.
Denote the resulting flow
$(A)(\VV',\Gamma,\mu)$
by $\WW=(f,w,\UU'')$.
We claim that this flow satisfies
(1)-(3) of 1.17. Indeed, $\WW$ is of radius $\mu$
and $\GG (\WW)\leq A$
by definition.
Since $\WW$ is subordinate to $\VV$ we have
$N(\WW)\leq N(\VV)$.
By 1.16 we have $\Vert w\Vert\leq\Gamma\Vert v\Vert =C$.
Further, let $\gamma$ be a $w$-trajectory.
The time which $\gamma$ can spend in
$\cup_{\crr} (U_p(r_p)\setminus
U_p(\mu))$
is not more than
$6\cdot N(\VV')\cdot r_p /C'\leq 6\cdot
2^n\cdot n \mu_0/C'\leq\beta/2$
and the time which $\gamma$ can spend in
$M\setminus \cup_{p\in S(f)}U_p(r_p)$ is not more than
$R(n)\cdot C'/C$ (since in this set we have
$w=C/C'\cdot v$), which does not exceed $\beta/2$.
\quad$\square$.
\enddemo

\document

\head{\S2.Transversality notions}\endhead
In this section we study the families of submanifolds of certain type.
The main example of such families is provided by the union
 of descending discs of an almost good gradient of a Morse function.
\subhead{A. Stratified submanifolds}\endsubhead
Let $\aaaa =\{A_0,...,A_k\}$
be a finite sequence of subsets of
a topological space $X$.
For $0\leq s\leq k$
we denote  $A_s$ also by $\aaaa_{(s)}$.
For $0\leq s\leq k$
we denote  $A_0\cup...\cup A_s$
by $A_{\leq s}$ and
also by $\aaaa_{(\leq s)}$.
We say that $\aaaa$ is a
\it compact family \rm
if  $\aaaa_{(\leq s)}$
is compact for every $s:0\leq s\leq k$.

\definition{Definition 2.1}
Let $M$ be a manifold without boundary.
 A finite sequence $\xx = \break \{X_0,...,X_k\} $
 of subsets of $M$ is called  \it
$s$-submanifold of $M$ \rm
($s$ for stratified) if
\roster
\item  $\xx$ is
disjoint and each
 $X_i$ is a submanifold of $M$ of dimension $i$
with the trivial normal bundle.
\footnote"*"{This restriction is technical, it makes proofs easier.}
\item $\xx$ is a compact family. \quad$\triangle$
\endroster
\enddefinition
For an $s$-submanifold $\xxx$ we denote $k$ by $\dim \xx$.
 For a diffeomorphism $\Phi :M\to N$ and an $s$-submanifold
 $\xx$ of $M$ we denote by $\Phi (\xx)$ the $s$-submanifold
of $N$ defined by $\Phi (\xx)_{(i)} = \Phi(\xx_{(i)})$.
\vskip0.1in
If $V$ is a submanifold of $M$ and $\xx$ is an $s$-submanifold of $M$,
then we say that $V$ \it is transversal to $\xx$ \rm
(notation: $V\pitchfork \xx$) if $V\pitchfork \xx_{(i)}$ for each $i$.
If $V$ is a compact submanifold of $M$, transversal to an
 $s$-submanifold $\xx$, then the family $\{\xx_{(i)}\cap V\}$
is an $s$-submanifold of $V$ (and of $M$) which will be denoted
 by $\xx \cap V$.
\vskip0.1in
If $\xx , \yy$ are two $s$-submanifolds of $M$, we say, that $\xx$
 is \it transversal to $\yy$ \rm
 (notation: $\xx \pitchfork \yy$) if
$\xx_{(i)}\pitchfork \yy_{(j)}$ for every $i,j$; we say that $\xx$
is \it almost transversal to \rm
   $\yy$ (notation: $\xx\nmid\yy$) if
$\xx_{(i)}\pitchfork \yy_{(j)}$ for $i+j<\dim M$.

\remark{Remark 2.2}
$\xx\nmid\yy$ if and only if $ \xx _{(\leq i)} \cap  \yy _{(\leq j)} =
 \emptyset$
whenever $i+j<\dim M$.\quad$\triangle$
\endremark

\proclaim{Lemma 2.3}
Let $f:M\to\RRR$ be a Morse function, where $M$ is a
 closed manifold,
$v$ be an almost good $f$-gradient.
Denote $\dim M$ by $n$.
Then the family
$\{K(\inde i;v)\}_{0\leq i\leq n}$
is an $s$-submanifold
transversal to $f^{-1}(\lambda)$
for every regular value $\lambda$ of $f$.
\endproclaim
\demo{Proof}
The set $\{K(\indl i;v)\}$
 is compact (by Lemma 1.6)
and the normal bundle to $D(p,v)$ is obviously trivial.
$\square$
\enddemo
The $s$-submanifold
 $\{K(\inde i;v)\}_{0\leq i\leq n}$
 will be denoted
by $\kk (v)$, the $s$-submanifold $\kk (v) \cap f^{-1}(\lambda)$
by $\kk _\lambda (v)$.
\vskip0.1in
For a manifold $M$ without
boundary we denote
by $\vem$ the subspace of
\break
$\text{Vect}^{\infty}(M\times [0,1])$,
consisting of all the
$\smo$
 vector fields which
have the second coordinate
 zero and which vanish with all the partial
derivatives in $M\times \{0,1\}$.
Assume that $M$ is compact.
Then there is a
natural topology
on the set $\text{Vect}^{\infty}(M\times [0,1])$ (see,
e.g. [19,\S 8]), with respect to which $\vem$ is a
closed subset. Further, for
every $v\in\vem$ and $t\in [0,1]$
the map $x\mapsto\gamma(x,t;v)$
is a $\smo$ diffeomorphism of $M$,
denoted by $\Phi(v,t)$.

The following transversality result
is proved by induction in
$\dim\xx +\dim\yy$
with the help of standard general
position argument like [6, Lem.5.3].
We omit the proof and just indicate
that the openness of the set in consideration follows from
Remark 2.2.

\proclaim{Theorem 2.4}
Let $\xx ,\yy$ be $s$-submanifolds
of a closed manifold $M$.
Then the set of $v\in\vem$, such that
$\Phi(v,1)(\xx)\nmid\yy$, is open and
dense in $\vem$.\quad$\square$
\endproclaim

\subhead{B. Good fundamental systems of neighborhoods and $ts$-submanifolds}\endsubhead
\definition{Definition 2.5}
Let $X$ be a topological space, $A\subset X$
 a closed subset, $I$ an open interval
$]0,\delta_0[$. A
\it
good fundamental system of neighborhoods of $A$
\rm
(abbreviation: $gfn$-system for $A$)
is a family
$\{ A(\delta )\}_{\delta\in I}$
of open subsets of $X$,
 satisfying the following conditions:
\roster
\item"(fs1)"  for each $\delta\in I$ we have $A\subset A(\delta)$ and
$\delta_1<\delta_2\Rightarrow A(\delta_1)\subset A(\delta_2)$,
\item"(fs2)"
for each $\delta\in I$ we have
 $\overline{A(\delta)}=\bigcap\limits_{\theta>\delta}A(\theta)$,
\item"(fs3)" $A=\bigcap\limits_{\theta>0}A(\theta)$.
\quad$\triangle$
\endroster\enddefinition
\proclaim{Lemma 2.6} Assume that $X$ is compact. Let
$\{ A(\delta )\}_{\delta\in I}$
be a $gfn$-system for $A$.
Then: \roster\item      The family
$\{ A(\theta )\}_{ \theta>0}$
is a fundamental system of neighborhoods of $A$.
\item $\forall\delta\in I$ the family
$\{ A(\theta )\}_{ \theta>\delta}$
is a fundamental system of neighborhoods of $\overline{A(\delta)}$.
\endroster
\endproclaim
\demo{Proof}
(1) Let $U$ be an open neighborhood of $A$.
The sets $\{X\setminus \overline{A(\theta)}\}_{\theta >0}$
form an open covering of the compact $X\setminus U$.
There is a finite subcovering, therefore there is $\theta >0$, such that
$\overline{A(\theta)}\subset U$.
(2) is proved similarly.\quad$\square$
\enddemo

\definition{Definition 2.7} Let $X$ be a topological space,
$\aaaa = \{A_0,...,A_k\}$ be a compact family
of subsets of $X$~,~ $I$  an open interval
$]0,\delta_0[$. A
\it
good fundamental system of neighborhoods of $\aaaa$
\rm
(abbreviation: $gfn$-system for $\aaaa$)
is a family
$\{ A_s(\delta )\}_{\delta\in I,0\leq s\leq k}$
of open subsets of $X$, satisfying the following conditions:
\roster
\item"(FS1)" For every $s : 0\leq s\leq k$ and
 every $\delta\in I$ we have $A_s\subset A_s(\delta)$
\item"(FS2)" For every $s : 0\leq s\leq k$
 we have $\delta_1<\delta_2\Rightarrow
A_s(\delta_1)\subset A_s(\delta_2)$
\item"(FS3)" For every $\delta\in I$  and every $j$
with $0\leq j\leq k$
we have
$\overline{A_{\leq j}(\delta)} =
{\bigcap}_{\theta >\delta}\big( A_{\leq j}(\theta)\big)$
\item"(FS4)" For every   $j$
with $0\leq j\leq k$ we have
$ A_{\leq j} =
{\bigcap}_{\theta >0}\big( A_{\leq i}(\theta)\big)$.
 \endroster
$I$ is called \it interval of definition
\rm
of the system.
\quad$\triangle$
\enddefinition

The following  lemma follow  from 2.6.
\proclaim{Lemma 2.8} Assume that $X$ is compact,
and let $\{ A_s(\delta )\}_{\delta\in I ,0\leq s\leq k}$
be a $gfn$-system for a compact family $\aaaa$.
Then for every $s : 0\leq s\leq k$
\roster
\item $\forall\delta\in I$ the family
$\{ A_{\leq s}(\theta )\}_{ \theta>\delta}$
is a fundamental system of neighborhoods of $\overline{A_{\leq s}(\delta)}$.
\item      the family
$\{ A_{\leq s}(\theta )\}_{ \theta>0}$
is a fundamental system of neighborhoods of $A_{\leq s}$.
\quad$\square$
\endroster
\endproclaim

Let $I=]0,\delta_0[$ and $Z=\{Z(\delta)\}_{\delta\in I}$
be a family of subsets of some space $X$.
Let $0<\epsilon_0<\delta_0$.
The family $\{Z(\delta)\}_{\delta\in ]0,\epsilon_0[}$
will be called \it restriction  of $Z$. \rm
We adopt the same terminology also for
good fundamental systems of neighborhoods
of compact families.

The basic example of a $gfn$-system is
 given by the next lemma, which follows
immediately from  Proposition 1.6.
\proclaim{Lemma 2.9} Let $M$ be a
 closed manifold or a compact
cobordism. Let $f:M\to\RRR$ be a Morse function
and $v$ be an almost good $f$-gradient.
Then for some $\epsilon>0$
 the family
$\{ B_\delta (\inde s;v)\}_{\delta\in ]0,\epsilon[,0\leq s\leq n}$
is a $gfn$-system for $\kk (v)$.
\quad$\square  $
\endproclaim

 The $gfn$-system,
introduced above,
 will be denoted
by $\KKK (v)$. Note that there is no canonical choice
of the interval  of definition for this system.
We shall say that $\lambda>0$ is in the interval of definition
of $\KKK (v)$ if there is
$\epsilon>\lambda$ such that
$\{ B_\delta (\inde s;v)\}_{\delta\in ]0,\epsilon[,0\leq s\leq n}$
is a $gfn$-system.
1.6 implies that if $v$ is
$\lambda$-separated, then
$\lambda$ is in the interval of definition
of $\KKK(v)$.
\vskip0.1in
                             Let $M$ be a manifold without
boundary, $\xx$ be an $s$-submanifold of $M$.
A good fundamental system of neighborhoods of $\xx$
will be called
\it thickened stratified
submanifold with the core $\xx$ \rm
(abbreviation: $ts$-submanifold).
For a $ts$-submanifold
$\XXX =\{ X_s(\delta )\}_{\delta\in I,0\leq s\leq k}$
we shall denote $X_i(\delta)$ by $\XXX _{(i)}(\delta)$,
and $X_{\leq i}(\delta)$
  by
$\XXX _{(\leq k)}(\delta)$.

 Lemma 2.9 implies that if $f:M\to\RRR$
is a Morse function on a closed manifold and $v$ is
an almost good $f$-gradient,
then $\KKK (v)$ is a $ts$-manifold with the core
$\kk (v)$.

\subhead{C. Tracks of subsets, $s$-submanifolds and $ts$-
submanifolds}\endsubhead
Let $f:W\to[a,b]$ be a Morse function on a compact riemannian
 cobordism
 $W$, $f^{-1} (b)=V_1,
f^{-1} (a)=V_0$, $v$ be an $f$-gradient.

 Let $X\subset V_1$.
 The set $\{\gamma (x,t;-v)\mid t\geq 0,x\in X\}$ will be
called \it track \rm of $X$ (with respect to $v$)
and denoted by $T(X,v)$.

\proclaim{Lemma 2.10}  \roster\runinitem If $X$ is compact,
 then  $T(X,v)\cup K(v)$ is compact.
\item  If $X$ is compact, and
every $(-v)$-trajectory starting from a point of $X$
reaches $V_0$ then $T(X,v)$ is compact.
\item For any $X$ we have $\overline{T(X,v)\cup K(v)} =
T(\overline X,v)\cup K(v)$.
\item  For any $X$
and $\delta>0$
 we have $\overline{T(X,v)\cup B_\delta (v)} =
T(\overline X,v)\cup D_\delta (v)$.
 \endroster
\endproclaim
\demo{Proof}
(1) The set $Y=W\setminus (T(X,v)\cup K(v))$
 consists of such
points $y\in W$ that $\gamma (y,t;v)$
 reaches $V_1$ and meets it at a point which
is not in $X$. Then
5.4 implies that $Y$ is open.

  (2) By an easy compactness
argument there is $\delta>0$ such that
$T(X,v)\cap B_\delta (p) =\emptyset$
for every $p\in S(f)$. Therefore
$T(X,v)\cap B_\delta (v)=\emptyset$ which
together with (1) implies (2).

(3) We have obviously:
$T(X,v)\cup K(v)\subset
 T(\overline X,v)\cup K(v)\subset
\overline{T(X,v)\cup K(v)}$, which implies (3) in view of (1).
(4) is proved similarly.\quad$\square$
\enddemo
 Note  that if $X$ is a submanifold of $V_1$
 of dimension $k$, then
$T(X,v)\cap\Wkr$ is a submanifold of
$\Wkr$ of dimension $k+1$.
If $\lambda$ is a regular
 value of $f$,
 then  $T(X,v)\cap f^{-1}(\lambda)$ is a submanifold
of dimension $k$ of $f^{-1}(\lambda)$.

From here to the end of the section we assume
that $v$ is an almost
 good $f$-gradient.
$\aaaa$ stands for an $s$-submanifold
$\{A_0,...,A_k\}$ of $V_1$, such that
$\aaaa\nmid\kk _b(-v)$.

 \definition{Definition  2.11}
For $0\leq i\leq k+1$ denote by $TA_i(v)$ the set
$T(A_{i-1} ,v)\cup K(\inde i ;v)$ (where we set
$A_{-1}=\emptyset$).
Denote by $\tttt (\aaaa ,v)$ the
 family
$\{ TA_i(v) \}_{0\leq i\leq k+1}$
 of subsets of $W$
and by $\stind vb\lambda (\aaaa)$ the family
$\{ \big(TA_{i+1}
(v) \big)
\cap f^{-1}(\lambda)\}_{0\leq i\leq k}$
of subsets of $f^{-1}(\lambda)$.
If the values $b,\lambda$ are clear from the context
 we shall abbreviate
$\stind vb\lambda (\aaaa)$ to $\st v (\aaaa)$.
The family $\tttt (\aaaa ,v)$ will be called \it track \rm
of $\aaaa$, and the family $\stind vb\lambda (\aaaa)$
will be called \it $\st v$-image \rm of $\aaaa$.
\quad$\triangle$
\enddefinition
\proclaim{Lemma 2.12} \roster\runinitem $\tttt (\aaaa ,v)$
and
$\stind vb\lambda (\aaaa)$ are compact families .
\item If $\lambda$ is a regular value of $f$, then
$\stind vb\lambda (\aaaa)$
is an $s$-submanifold of  $f^{-1}(\lambda)$.
\endroster\endproclaim
\demo{Proof}
(1) Let
$0\leq s\leq k+1$. Denote
$ T(A_{\leq s-1} ,v)\cup K(\indl s;v)$
by $Y(s)$.
Let $\phi :W\to[a,b]$ be an ordered Morse
function, adjusted to $(f,v)$.
Let $\lambda$ be a regular value of
 $\phi$, such that all the critical points
 of $\phi$ of indices $>s$ are
 above $\lambda$ and all the critical
points of $\phi$ of indices $\leq s$
are below $\lambda$.
 Since
$\aaaa\nmid (\kk (-v)\cap V_1)$, every
 $(-v)$-trajectory, starting at a point
of $A_{\leq s-1}$ reaches
$\phi^{-1}(\lambda)$. Then
$Y(s)\cap\phi^{-1}([\lambda , b])$
is compact by 2.10(2)
and
$Y(s)\cap\phi^{-1}([a, \lambda ])$
is compact by 2.10(1), therefore
$Y(s)$ is compact.
(2) is obvious.\quad$\square$
\enddemo

We shall now define the notion of  the track of
a $ts$-submanifold.

\definition{Definition 2.13} Let
 $\AAA =\{ A_s(\delta )\}_{\delta\in ]0,\delta_0[ ,0\leq s\leq k}$
be a $ts$-submanifold of $V_1$ with the core
$\aaaa$.
Assume that $v$ is $\delta_1$-separated (where $\delta_1>0$).
 For $0<\delta<\min (\delta_0,\delta_1)$
and $0\leq s\leq k+1$
set
$TA_s(\delta , v) = T(A_{s-1} (\delta),v)
\cup
 B_\delta (\inde s; v)$ (where we set by definition
$A_{-1}(\delta)=\emptyset$).
\quad$\triangle$
\enddefinition
We shall prove  that some restriction of
$\{TA_s(\delta,v)\}$ is a $gfn$-system for
$\tttt (\aaaa ,v)$.
Up to the end of this section
$Q(s,\delta)$ stands for
$T(\overline{A_{\leq s-1}(\delta)},v)
\cup  D_\delta (\indl s;v)$.

\proclaim{Lemma 2.14}
Let $0<\epsilon\leq\min(\delta_0,\delta_1)$.
Then (i)$\Leftrightarrow$(ii).
\roster\item"(i)" For every $0<\delta<\epsilon$ and every
$s:0\leq s\leq k$ the set
$Q(s,\delta)$ is compact
\item"(ii)"
$TA_s(\delta , v)_{0<\delta<\epsilon,
0\leq s\leq k+1}$
is a $gfn$-system for $\tttt (\aaaa,v)$
\endroster\endproclaim
\demo{Proof} It is easy to see that for every
$\delta\in ]0,\min(\delta_0,\delta_1)[$ we have:
$$\gather
Q(s,\delta)={\cap}_{\theta>\delta} TA_{\leq s}(\theta,v),
\qquad
TA_{\leq s}(v)={\cap}_{\theta>0} TA_{\leq s}(\theta,v)
\tag2.1\\
TA_{\leq s}(\delta,v)\subset Q(s,\delta)
\subset \overline{TA_{\leq s}(\delta,v)}
\tag2.2
\endgather
$$
Now to prove (ii)$\Rightarrow$(i) note that
(FS3) together with (2.1)  implies
$Q(s,\delta)=\overline{TA_{\leq s}(\delta,v)}$.
To prove (i)$\Rightarrow$(ii) note that
(FS1), (FS2), and (FS4) hold for every
$\delta\in ]0,\min(\delta_0,\delta_1)[$,
and (i) implies (FS3) in view of
2.2.\quad$\square$\enddemo

\proclaim{Proposition 2.15}
There is $\epsilon \in ]0,\min (\delta_0,\delta_1)[$
such that
$\{TA_s(\delta , v) \}_{\delta\in]0,\epsilon[,~0\leq s\leq k}$
is a gfn-system for $\tttt (\aaaa ,v)$.
\endproclaim
\demo{Proof}
Let $\phi:W\to[a,b]$ be the Morse function
with respect to which $v$
is $\delta_1$-separated, and let
$\mu$
 be a regular
value of $\phi$ such that for every
 critical point $p$ of $\phi$ of indices $\leq s$
(resp. $>s$) the disc $D_p(\delta_1)$
is in $\phi^{-1}([a,\mu[)$
(resp. in
$\phi^{-1}(]\mu,b])$).
Since $\aaaa\nmid\kk  _b(-v)$, every
$(-v)$-trajectory
starting at a point of $A_{\leq s-1}$
reaches $\phi^{-1}( \mu)$.
$A_{\leq s-1}$ is compact and
$\{A_{\leq s-1}(\theta)\}_{\theta>0}$
is a $gfn$-system for $A_{\leq s-1}$;
therefore an easy compactness argument
based on 5.4 shows that there is $\theta>0$
such that
each $(-v)$-trajectory
starting at a point of
$A_{\leq s-1}(\theta)$
reaches $\phi^{-1}(\mu)$.

We claim that $\epsilon=\min(\delta_0,\delta_1,\theta)$
satisfy the conclusions of our proposition.
Indeed, let $\delta<\epsilon$.
Then
each $(-v)$-trajectory
starting at a point of
$\overline{A_{\leq s-1}(\delta)}$
reaches $\phi^{-1}(\mu)$, and
$T(\overline{A_{\leq s-1}(\delta)},v)$
is compact by 2.10(2).
Denote
$\stind vb\mu \big(\overline{A_{\leq s-1}(\delta)}\big)$
by $C(\delta)$; it is a compact set.
The intersection of the set
$Q(s,\delta)$
with $\phi^{-1}([\mu,b])$
is compact by the above, and its intersection with
$\phi^{-1}([a,\mu])$
is compact by 2.10(4), applied to the cobordism
$\phi^{-1}([a,\mu])$
and the compact $C(\delta)\subset\phi^{-1}(\mu)$.
\quad$\square$\enddemo

 We shall denote the $gfn$-system,
introduced in 2.15, by $\TTT(\AAA,v)$
and call it \it track \rm of $\AAA$.
Note that there is no canonical choice of
the interval of definition
for this system. We shall say that
$\lambda>0$ belongs to the interval of definition
of  $\TTT(\AAA,v)$
if there is
$\epsilon>\lambda$ such that
\break
$\{TA_s(\delta , v) \}_{\delta\in]0,\epsilon[,~0\leq s\leq k}$
is a $gfn$-system.
The next two lemmas contain
some properties of $\TTT(\AAA,v)$.

\proclaim{Lemma 2.16}
Assume that $\delta>0$ is in the interval of definition
of $\TTT(\AAA,v)$. Then for $\lambda=a,b$ we have:
$\overline{\TTT(\AAA)_{(\leq s)}(\delta)}
\cap f^{-1}(\lambda)=
\overline{\TTT(\AAA)_{(\leq s)}(\delta)
\cap f^{-1}(\lambda)}$.
\endproclaim
\demo{Proof} Obvious. \quad$\square$\enddemo

Let $\alpha,\beta$ be regular values of $f$,
such that $a<\alpha<\beta<b$ and
there are no critical points of $f$ in
$f^{-1}([a,\alpha]\cup[\beta,b])$.
Denote by $\widetilde f$, $\widetilde v$
the restrictions of $f$,$v$ to
$\widetilde W =f^{-1}([\alpha,\beta])$.
Denote $f^{-1}(\alpha)$ by $\widetilde V_0$,
$f^{-1}(\beta)$ by $\widetilde V_1$.
Let
$\widetilde\AAA=\{\widetilde A_s(\delta)\}_
{\delta\in]0,\delta_0[, 0\leq s\leq k}$
be a $ts$-submanifold of $\widetilde V_1$ with the core
$\widetilde\aaa=\stind vb\beta (\aaa)$, and
assume that for every $s:0\leq s\leq k$ and
for every $ \delta\in ]0,\delta_0[$
we have:
$\stind vb\beta (A_s(\delta))\subset \widetilde A_s(\delta)$.
Assume that $\widetilde v$ is
$\delta_1$-separated.

\proclaim{Lemma 2.17}
Let $0<\nu<\min(\delta_0,\delta_1)$
and assume that $\nu$ is in the interval of definition
of $\TTT(\widetilde\AAA
,\widetilde v)$. Then $\nu$ is in the
interval of definition of $\TTT(\AAA,v)$.
\endproclaim
\demo{Proof}
Choose $\theta$ in $]\lambda,\min(\delta_0,\delta_1)[$.
Since $v$ is obviously $\delta_1$-separated,
it is sufficient to prove that for every
$\delta\in]0,\theta[,~s:0\leq s\leq k$
the set
$Q(s, \delta)$
is compact.
The set $Q(s, \delta)\cap f^{-1}([\beta,b])$ is compact
by 2.10(2).
The set $Q(s, \delta)\cap f^{-1}([\alpha, \beta])$ equals to
$T\big(\overline{\stind vb\beta (A_{\leq s-1}(\delta))},
\widetilde v\big)
\cup
D_\delta(\indl s;\widetilde v)$.
If $Z'\subset Z$ are compacts in $\widetilde V_1$ then
$T(Z',\widetilde v)$
is closed in $T(Z,\widetilde v)$, therefore
 $Q(s, \delta)\cap f^{-1}([\alpha, \beta])$
is a closed subset of
$T(\overline{\widetilde A_{\leq s-1}(\delta)},\widetilde v)
\cup
D_\delta(\indl s;\widetilde v)$,
therefore it is compact.
Finally,
    compactness of
$Q(s, \delta)\cap f^{-1}([a,\alpha])$
follows from
the compactness of
$Q(s, \delta)\cap f^{-1}([\alpha,\beta])$
and from the absence of critical points
of $f$ in $f^{-1}([a,\alpha])$.\quad$\square$
\enddemo

\document

\head{\S3. S-Construction}\endhead
In Subsection A we present a construction, which
produces from a Morse function without critical points
on a compact cobordism
another one,
 having two series of critical points
and behaving roughly as "skladka". This new
function is equipped with an almost good
gradient, and the main property of the
 construction is an explicit
estimate of the quickness of the resulting flow.
The properties of
S-construction are
listed in the theorem 3.1 below.
We invite the reader first to have a look
at the construction
itself (proof of 3.1). In Subsection B we
prove
 $1.17(n)\Rightarrow 1.18(n+1)$
with the help of the S-construction.

\it Terminology: \rm
If $W$
is a compact cobordism, $v\in \text{Vect}^1(W,\bot)$
and $U\subset W$,
we say that $U$ is
\it $v$-invariant \rm
if for every $x\in U$
the trajectory $\gamma(x,\cdot ;v)$ is defined on
$[0,\infty [$ and
$\gamma(x,t;v)\in U$ for all $t\geq 0$. We say that
$U$ is
\it weakly $v$-invariant \rm
if for every $x\in U$
we have: $\gamma(x,t;v)\in U$
for all $t$ of the interval of definition of $\gamma(x,\cdot;v)$.

\subhead A. S-construction \endsubhead
Let $W$ be a compact riemannian cobordism.
During \S 3 we denote by $\vert\cdot\vert$
the norm on the tangent spaces, induced
by the metric and by $\Vert\cdot\Vert$
the corresponding $C^0$ norm
in the space of vector fields.

\proclaim{Theorem 3.1}
Let $g:W\to [a,b]$ be a Morse function
without critical points,
$g^{-1}(a) =V_0, g^{-1}(b) =V_1$.
Let $C>0$ and let $w$ be a $g$-gradient,
such that $\nrw \leq C$.
Denote $g^{-1}(\frac {2a+b}3)$ by $V_{1/3}$,~
$g^{-1}(\frac {2b+a}3)$ by $V_{2/3}$,~
$g^{-1}(\frac {a+b}2)$ by $V_{1/2}$.
Denote
$g^{-1}([a,\frac {2a+b}3])$ by $W_0$,
$g^{-1}([\frac {2b+a}3,b])$ by $W_2$,
$g^{-1}([\frac {2a+b}3,\frac {2b+a}3])$
by $W_1$.
Denote by $\text{grad}(g)$ the riemannian
gradient of $g$ and by $\text{grd}(g)$ the
 vector field
$\text{grad}(g)/
\vert \text{grad}(g) \vert$.
Denote by $T$ the maximal length
of the domain of a $\text{grd}(g)$-trajectory.
Then there is $\nu_0>0$ such that:

For every $s$-submanifold $\xx$ of $V_0$,
$s$-submanifold $\yy$ of $V_1$,
every almost good $\MM$-flow
$\VV_1=(F_1, u_1, \UU_1)$
on $\VODIN$ and
almost good
$\MM$-flow $\VV_2=(F_2, u_2, \UU_2)$
on $\VDVA$, and every
$\mu\leq \min (\nu_0,d(\VV_1),d(\VV_2))$
there is an almost good
$\MM$-flow $\VV=(F,u,\UU)$ on $W$,
satisfying the following properties:
\roster\item
$F:W\to[a,b],\quad V_0=F^{-1}(a),
V_1=F^{-1}(b), \quad S(F)=S(F_1)\cup S(F_2)$;
in a neighborhood of $\partial W$
we have: $F=g,u=w$.

\item    $\VV $ is of radius $\mu$;\quad
      $\GG (\VV)\leq {\frac {3}{2}}
\max (\GG (\VV _1), \GG (\VV _2))$.

\item $\VODIN , \VDVA, W_1 $ are
$(\pm u)$-invariant.
$W_0$
 is $u$-invariant and weakly $(-u)$-invariant.
 $W_2$ is $(-u)$-invariant and weakly
$u$-invariant.
\item
Assume that $\VV _1$ is $(C_1,\beta_1,\mu/2)$-quick
and  $\VV _2$ is $(C_2,\beta_2,\mu/2)$-quick.
Then $\VV$ is
$(3/2(C+C_1 +C_2),~\beta_1 +
\beta_2+5+4T/C,~\mu)$-quick.
\item
Let $\gamma$ be an $u$-trajectory.
If $     \Im \gamma\subset W_1$, then
$N(\gamma,\mu)\leq N(\VV _1,\mu)+N(\VV _2,\mu)$.
If $     \Im \gamma\subset W_0$, then
 $N(\gamma,\mu)\leq N(\VV _1,\mu)$.
If $     \Im \gamma\subset W_2$, then
$N(\gamma,\mu)\leq N(\VV _2,\mu)$.
\item
$\kk_a(u)\nmid \xx;\quad \kk_b(-u)\nmid\yy.$
\endroster
\endproclaim
\subsubhead{Proof}\endsubsubhead
The proof occupy the rest of Subsection A.
We shall assume that $a=0, b=1$, since the
general case is easily reduced to this one
by an affine transformation of $\RRR$.
\subsubhead{1. Function $f$, its gradient $v$,
and the choice of $\mu_0$}\endsubsubhead

\proclaim{Lemma 3.2}
There is a Morse function $f:W\to [0,1]$
 without critical points and an $f$-gradient
$v$, such that:
\roster\runinitem  $\nrv \leq C$;
\item for $x$ in a neighborhood of
$V_0\cup V_1$ we have: $ f(x)=g(x)$
and $ v(x)=w(x)$;
\item $f^{-1}(1/3)=V_{1/3},~ f^{-1}(2/3)=V_{2/3},~
f^{-1}(1/2)=V_{1/2}$;
\item for $x$ in a neighborhood of
 $V_{1/3}\cup  V_{1/2}\cup V_{2/3}$
we have:
\break
$ \vert v(x)\vert =C\aand df(v)(x)=C$;
\item for
$\lambda=1/3,~1/2,~2/3\aand x\in g^{-1}(\lambda)$
we have: $v(x)\bot g^{-1} (\lambda)$;
\item The maximal length of the domain of a
 $v$-trajectory is not more than $2T/C$.
\endroster
\endproclaim
\demo{Proof}
Let $U$ be an open neighborhood of
 $V_0\cup V_1$, such that
$U\cap(\VODIN\cup\VDVA\cup\VM)=\emptyset$
and let $h:W\to [0,1]$ be a $\smo$ function such
that $\supp ~h\subset U$ and for $x$ in a
neighborhood of $\partial W$ we have $h(x)=1$. Set
$v(x)=h(x)w(x) + (1-h(x))C\text{grd}(g)$. It is
 obvious that $v$ is a $g$-gradient, satisfying
 1) and 5).
We have also $v(x)=w(x)$ nearby $\partial W$,
as well as $\nrv = C$ nearby
$\VODIN\cup\VDVA\cup\VM$.
6) holds also, if only $U$ was chosen sufficiently small.
Applying now Corollary  8.14 of [19] to the cobordisms
$W_0,~ g^{-1}([1/3,1/2]),
g^{-1}([1/2,2/3]),~ W_2$  and the restrictions of $g$ and $v$
to these cobordisms and glueing the results together
we obtain a Morse function $f:W\to [0,1]$, satisfying
(together with $v$) all of our conclusions.\quad$\square$
\enddemo

For $\lambda\in[0,1]$ we denote $f^{-1}(\lambda)$ by
 $V_\lambda$. Fix some
$\epsilon \in ]0,1/12[$.
For $\nu>0$ sufficiently small
the map
$(x,\tau)\mapsto\gamma(x,\tau;v/C)$ is defined on
$V_i\times[-2\nu,2\nu]$, where $i=1/3, 1/2, 2/3$, on
$V_0\times[0, 2\nu]$
and on $V_1\times[-2\nu, 0]$. The corresponding
embeddings will be denoted by
$$\gather
\Psi_0(\nu) : V_0 \times [0,2\nu ] \to W,\quad
\Psi_{3/2}(\nu) : V_{1/2} \times [-2\nu,2\nu ] \to W,\quad
 \Psi_3(\nu) : V_1 \times [-2\nu,0] \to W ,\\
\Psi_1(\nu) : V_{1/3} \times [-2\nu, 2\nu] \to W,  \qquad
\Psi_2(\nu) : V_{2/3} \times [-2\nu, 2\nu ] \to W.
\endgather
$$
\it Terminology: \rm
Two riemannian metrics on a manifold
$N$ are called
\it $C$-equivalent, \rm  if
for every tangent vector $h$ to $M$
we have:
$\vert h\vert_{g_1}/\vert h\vert_{g_2} \leq C,\quad
\vert h\vert_{g_2}/\vert h\vert_{g_1} \leq C$.

Let $\nu_0$ satisfy the following restriction:

($\frak R$ ): $2\nu_0<C$.
For $i=0, 1/3, 1/2, 2/3, 1$ the riemannian metric induced
by $\Psi_{3i}(\nu_0)$ from $W$
and the product metric on the domain
of  $\Psi_{3i}(\nu_0)$ are $3/2$-equivalent;
further, we have:
$f(\Im\Psi_{3i}(\nu_0))\subset
 ]i -\epsilon, i +\epsilon[$.
For  $i=1/3,~1/2,~2/3$
we have:
$\vert v(x) \vert =df(v(x))=C$
for $x$ in a neighborhood of
$\Im\Psi_{3i}(\nu_0)$.

We shall prove that this $\nu_0$ satisfy
the conclusions of our theorem.
So let $\xx, \yy, \VV_1,
\VV_2$ be as in the
statement of the theorem. Denote
$\dim W$ by $n+1$. Set
$\UU_1=\{\Phi_{1p}:U_{1p}\to B^n(0,r_{1p})\}$ and
$\UU_2=\{\Phi_{2p}:U_{2p}\to B^n(0,r_{2p})\}$.
Let $\mu\leq\min(\nu_0, d(\VV_1), d(\VV_2))$.
We shall denote $\Psi_{3i}(\mu)$ by
$\Psi_{3i}$.

For $\theta\in]0,2\mu]$ and
 $i=1,3/2,2$ we denote by $\text{Tb}_i(\theta)$ the set
$\Psi_i ( V_{i/3}\times [-\theta,\theta])$.
For $\theta\in]0,2\mu]$
we denote
by $\text{Tb}_0(\theta)$ the set
$\Psi_0 ( V_0\times [0,\theta])$ and by
$\text{Tb}_3(\theta)$ the set
$\Psi_1 ( V_1\times [-\theta,0])$.

 For $\lambda\in [0,1/3]$ we denote by
 $\llll_\lambda (u_1)$ the
$s$-submanifold
$\stind v{1/3}\lambda (\kk (u_1))$ of $V_\lambda$.
For $\lambda\in [1/3,1]$ we denote
 by $\llll_\lambda (u_1)$
 the $s$-submanifold
 $\stind {(-v)}{1/3}\lambda \kk (u_1)$ of
$V_\lambda$.
For $\lambda\in [2/3,1]$ we denote
 by $\llll_\lambda (-u_2)$
 the $s$-submanifold $\stind {(-v)}{2/3}\lambda (\kk (-u_2))$
and
for $\lambda\in [0,2/3]$ we denote by
$\llll_\lambda (-u_2)$ the
 $s$-submanifold $\stind v{2/3}\lambda (\kk (-u_2))$.
Note that since $v$ have no zeros,
$\stind v\alpha\beta$ is a diffeomorphism
of $V_\alpha$ onto $V_\beta$
for any $\beta<\alpha$.
Note also
that $\llll_{2/3} (-u_2)$ equals to
$\kk (-u_2)$ and
$\llll_{1/3} (u_1)$ equals to
$\kk (u_1)$.

\subsubhead{2. Auxiliary functions and vector fields}
\endsubsubhead

Let $\chi:[0,1]\to [0,1]$ be a $\smo$
 function, with the following properties: $$\gather
\chi (x)=x\ffor x\in [0,\epsilon]\cup [1-\epsilon , 1];\qquad
\chi (x) =1-x\ffor x\in [1/2-\epsilon, 1/2+\epsilon];\\
\chi (x) = 2/3-(x-1/3)^2
\for
 x\in [1/3-\epsilon,1/3+\epsilon];\quad
\chi'(x)>0
\for
 x\in[0,1/3[~\cup ~]2/3,1];
\\
\chi (x) = 1/3+(x-2/3)^2\for
 x\in [2/3-\epsilon,2/3+\epsilon];\quad
\chi'(x)<0 \for x\in ]1/3 , 2/3[.
\endgather
$$
\proclaim{Lemma 3.3}
\roster\runinitem $(\Psi_i^{-1})_*(v)=(0,C)\ffor i=0,~1,~2,~3/2,~3$.
\item $(f\circ \Psi_i)(x,\tau) =\tau +i/3 ,
  \where i=  1,2,3/2  $.
\item $(\chi\circ f\circ\Psi_1)(x,\tau)=2/3-\tau^2,\quad
(\chi\circ f\circ\Psi_2)(x,\tau)=1/3+\tau^2$.
\item For $\lambda\in [-2\mu,2\mu]$ and
$i=1,3/2,2$ we have
$f^{-1}(i/3+\lambda)=
\Psi_i(V_{i/3}\times\{\lambda\}).$\endroster
\endproclaim
\demo{Proof} (1) - (3) follow immediately from ($\frak R$).
To prove (4) let $i=1$ and $\lambda\in [-2\mu,2\mu]$.
Each $v$-trajectory intersects each level surface
of $f$  at exactly one point.
Therefore on each $v$-trajectory starting at a point
$x\in V_{1/3}$(and, hence, on every
$v$-trajectory), the point
$\Psi_1(x,\lambda)$ is the only
point in $f^{-1}(1/3+\lambda)$.
The cases $i=3/2, 2$ are similar.
\quad$\square$
\enddemo

Let $B:\RRR\to [0,1]$ be a $\smo$ function
such that
$B(t)=0\ffor \vert t\vert\geq 5\mu/3\aand
B(t)=1\ffor  \vert t\vert\leq 4\mu/3.$
Let $B_1:\RRR\to\RRR^+$ be a $\smo$ function
such that $\supp B_1\subset ]\mu, 2\mu[$
and that $\int_0^{\infty} B_1(t)dt=C$.

Let $z_0$ be a $\smo$ vector field on $V_0$ such
that $\Phi(z_0,1)(\llll _0 (u_1))\nmid\xx$.
Let $z_1$ be a $\smo$ vector field on $V_1$ such
that $\Phi(z_1,1)(\llll _1 (-u_2))\nmid\yy$.
Let $z_{1/2}$ be a $\smo$ vector field on $V_{1/2}$
such that $\Phi(z_{1/2},1)(\llll _{1/2} (u_1))
\nmid\llll _{1/2} (-u_2)$.
We shall assume that $z_i$ (where $i=0,1/2,1$)
are chosen so small, that
$\sup_\tau\vert B_1(\tau)\vert
\cdot\Vert z_i\Vert<C/9$.

\subsubhead{3. Morse function $F$}\endsubsubhead
$$\gather
\text{\rm Let } a_1,a_2>0. \text{\rm ~~
 Set:}\quad
 F(y)
 =\psi (y) \ffor y\in W
\setminus (\Tb _1(2\mu)\cup\Tb _2 (2\mu));\\
(F\circ\Psi _1)(x,\tau)
  =
a_1 B(\tau)F_1(x) +2/3-\tau^2 \ffor (x,\tau)\in
 \VODIN\times [-2\mu,2\mu]; \\
(F\circ\Psi _2)(x,\tau)
  = a_2 B(\tau)F_2(x) +1/3+
\tau^2\ffor (x,\tau)\in \VDVA\times [-2\mu,2\mu]
\endgather
$$
Since $(\psi\circ\Psi_1)(x,\tau)=2/3 -\tau^2$ and
$(\psi\circ\Psi_2)(x,\tau)=1/3 +\tau^2$, these
formulas define correctly a smooth function
$F:W\to\RRR$, which equals to
$f$ nearby $\partial W$.
To find critical points of $F$ note that $S(F)$ is
contained obviously in
$\Tb _1(2\mu )\cup\Tb _2(2\mu)$. For
$(x,\tau)\in\VODIN\times[-2\mu,2\mu]$ we have
$d(F\circ\Psi_1)(x,\tau)=
(a_1 B'(\tau) F_1(x) -2\tau)d\tau +
a_1 B(\tau) dF_1(x)$.
For $a_1$ small enough this can
vanish
only for $\tau =0$. We conclude therefore
(applying the same reasoning to $\Tb _2(2\mu)$),
 that $S(F)=S(F_1)\cup S(F_2)$ if only $a_i$ are
small enough (and we make this assumption from now on).
To prove that $F$ is a Morse function we shall
explicit the standard charts for $F$.
 Let $p\in S(F_1)$ and write
$(F_1\circ\Phi_{1p}^{-1})(x)=
F_1(p)+
\sum_i\alpha_ix_i^2$.
 Consider the chart
$\Phi_{1p}\times\id :U_{1p}\times ]-\mu,\mu[\to B^n(0,r_{1p})\times
 ]-\mu,\mu[$ of the manifold $\VODIN\times
]-\mu,\mu[$ around the point $p\times 0$.
We have
$F\circ\Psi_1\circ (\Phi_{1p}\times\id)^{-1}(x,\tau) =
F_1(p)+2/3+a_1 \sum_i \alpha_ix_i^2 -\tau^2$,
 therefore the chart
$\big( (\Phi_{1p}\times\id)
\mid (\Phi_{1p}\times\id)^{-1} (B^{n+1}(0,\mu)\big)\circ
(\Psi_1\mid\text{Im}\Psi_1)^{-1}$
is a standard chart of radius $\mu$
for $F$ at $p\in S(F_1)$.
 These charts
together with the similar ones for $q\in S(F_2)$
give an $F$-chart-system of radius $\mu$.
We shall denote this system
by $\UU$.
Note also that $\ind _F p=\ind _{F_1}p+1$ for
$p\in S(F_1)$ and that $\ind _F q=\ind _{F_2}q$
 for $q\in S(F_2)$.
 Note that if $a_1$ and $a_2$ are chosen
sufficiently small, then
$F(\Tb _1(2\mu))\subset
[2/3 -\epsilon, 2/3 +\epsilon]$, and
$F(\Tb _2(2\mu))\subset [1/3 -\epsilon, 1/3 +\epsilon]$.
In this case also $F^{-1}(1/2) =V_{1/2} \cup
 V_{\theta} \cup V_{\theta'}$, where
$\epsilon <\theta <1/3 -\epsilon,\quad
2/3 +\epsilon <\theta' < 1-\epsilon$.
\subsubhead
{    4.  $F$-gradient $u$   }
\endsubsubhead

Define a vector field $u$ on $W$ as follows:
$$\gather
u(y)=v(y) \ffor y\in \big( W_0\setminus
(\Tb _1 (2\mu)\cup\Tb _0 (2\mu))\big) \cup
 \big(W_2\setminus (\Tb _2(2\mu)\cup
 \Tb_3(2\mu)\big);\\
u(y)=-v(y)\ffor  y\in  W_1\setminus
\big(\Tb _1 (2\mu)\cup\Tb _2 (2\mu) \cup
 \Tb _{3/2}(2\mu)\big);\\
\big( (\Psi _0^{-1})_* (u) \big) (x,\tau)=
(-B_1(\tau)z_0(x),~ C);\quad
\big( (\Psi _3^{-1})_* (u) \big) (x,\tau)=
(B_1(-\tau)z_1(x),~C);
\\
\big( (\Psi _1^{-1})_* (u) \big) (x,\tau)=
\big(B(\tau)u_1(x), ~-B(\tau)\tau -
(1-B(\tau))\cdot C\cdot\text{sgn} \tau\big);\\
\big( (\Psi _{3/2}^{-1})_* (u) \big) (x,\tau)=
(-B_1(-\tau)z_{1/2}(x),~
 -C);\\
\big( (\Psi _2^{-1})_* (u) \big) (x,\tau)=
 \big(B(\tau)u_2(x), ~B(\tau)\tau +
(1-B(\tau))\cdot C\cdot\text{sgn} \tau\big).\\
\endgather
$$

An easy computation
 using the definition of $F$ and $u$
shows that  these formulas
define correctly a
$\smo$ vector field on $W$, and that
$(F,u,\UU)$ is an $\MM$-flow on $W$
of radius $\mu$, if only $z_0, z_1, z_{1/2}$
are small enough (which assumption we
make from now on). It is also easy to see
that
$\supp (u-v)\subset
\big(\cup_{i=0,~1,~3/2,~2,~3}\Tb_i(2\mu)\big)$,
and that  $u$ is an $f$-gradient in
 $W_0\setminus V_{1/3}$
and $W_2\setminus V_{2/3}$; ~
$(-u)$ is an $f$-gradient in
$W_1\setminus (V_{1/3}     \cup V_{2/3}  )   $.

We claim that $(F, u, \UU)$
satisfies all the conclusions of 3.1.
(1) and (2) follow immediately
from the construction. To prove
(3) note that $V_{1/3}$ is a closed submanifold
of $\kr {W}$
 and $u$ is tangent to $\VODIN$;
therefore $\VODIN$ is $(\pm u)$-invariant
 (same for $\VDVA$).
 To prove that $W_1$ is $u$-invariant,
let $x\in\kr {W}_1$ and assume that
for some $t$ we have
$\gamma (x,t;u)\notin \kr {W}_1$.
The trajectory $\gamma(x,\cdot ;u)$
does not intersect $\VODIN\cup\VDVA$.
Consider a continuous function
$\phi :t\mapsto f(\gamma(x,t;u))$.
The domain of definition of $\phi$
is an interval (finite or infinite)
of $\RRR$, ~ $\phi$ never takes values
1/3, 2/3, and $\phi_0\in ]1/3,2/3[$. Therefore
$\Im \phi\subset]1/3,2/3[$
and $W_1$ is $u$-invariant.
The other assertions of  (3)
are proved similarly.

\subsubhead{ 5. Estimate of the quickness of $\VV$}
\endsubsubhead

To obtain the estimate  $\Vert u\Vert\leq C$
note that the inequality
$\vert u(x)\vert\leq
3/2(C+C_1+C_2)$ is to be checked
only for $x\in\cup_i\Tb_i(2\mu)$,
where it is a matter of a simple computation;
we leave it to the reader.
 To estimate the
time, which an $u$-trajectory spends outside
 $\UU (\mu)$, note first that for a trajectory, starting at
 a point of $\VODIN$ (resp. $\VDVA$),
this time is not more than $\beta _1$ (resp.
$\beta_2$),
since it is actually a $u_1$- (resp. $u_2$)-trajectory.
Now let  $\gamma$
be an
$u$-trajectory, passing by
a point of $\kr {W}_1$. By (4) and (5)
it stays in $\kr {W}_1$ forever.
 Since $u$ is a $(-f)$-gradient in $\kr {W}_1$,
the function $t\mapsto f(\gamma (t))$
is strictly decreasing and
$\lim_{t\to -\infty} f(\gamma (t)) = 2/3,\quad
\lim_{t\to \infty} f(\gamma (t)) =1/3$.
 This implies
 $\lim_{t\to -\infty} \gamma (t) =q\in S(F_2)$ and
$\lim_{t\to \infty} \gamma (t) =p\in S(F_1)$.
We shall estimate the time which $\gamma$
spends outside $\UU (\mu)$
between the various level surfaces of $f$.

1) $f(\gamma (t))\in[2/3 -\mu/2,2/3]$. \quad
By Lemma
3.3(4) this condition is equivalent to:
$\gamma (t)\in\Psi_2(V_{2/3}\times[-\mu/2,0[)$.
The curve $\Psi_2^{-1}(\gamma(t))$
is a product of
an $u_2$-trajectory and the curve
$\tau\mapsto \alpha e^{\tau}$
with $\alpha<0$.
 Since
 $\Psi_2^{-1}(\UU(\mu))$
contains
$U_p(\mu/2)\times ]-\mu/2,\mu/2[$
for every $p\in S(F_2)$
the time which $\gamma (t)$ spends in
$\Psi_2(V_{2/3}\times [-\mu/2,0[)$
outside $\Psi_2^{-1}(\UU(\mu))$
is not more than $\beta_2$.

2) $f(\gamma (t))\in[2/3 -2\mu,2/3-\mu/2]$, in
other words,
$\gamma(t)\in \Psi_2(V_{2/3}\times [-2\mu,-\mu/2])$.
The vector field $(\Psi_2^{-1})_*(u)$
equals to
$(B(\tau) u_2(x),\varkappa (\tau))$, where
$\varkappa (\tau)\leq\tau$. Therefore the total time
which $\gamma$ can spend in this domain
is not more than
$\ln \big(2\mu/(\mu/2)\big) < 2$.

3) $f(\gamma (t))\in[1/2, 2/3-2\mu]$.\quad
In the domain
 $f^{-1}([1/2,2/3-2\mu])$
we have $u=v$. Therefore the time is
$\leq 2T/C$.

4) $f(\gamma (t))\in[1/2-2\mu, 1/2]$.\quad
Here $\gamma (t)\in\Psi_{3/2}(V_{1/2}\times
[-2\mu, 0])$. The second coordinate
of $(\Psi_{3/2}^{-1})_*(u)$ is equal to $-C$,
therefore the total time
spent here is not more than $2\mu/C\leq 1$.

Similarly to the cases 1) - 3) above
one shows that the time which
$\gamma$ spends in  \break
$f^{-1}(]1/3,1/2 -2\mu])\setminus \UU (\mu)$
is not more than
$2T/C+2+\beta_1$. Summing
 up, we obtain that the time
which $\gamma$
can spend in $\kr {W}_1\setminus\UU (\mu)$
is not more than
$\beta_1+\beta_2+4T/C+5$.

Similar analysis of behaviour of
 $u$-trajectories in $W_0$ and $W_2$
shows that the same estimate
holds in these cases also.

\subsubhead{6. Estimate of $N(\gamma,\mu)$}
\endsubsubhead

Let $\gamma(\cdot)$ be an $u$-trajectory in $W_1$.
If $\gamma$ is in $V_{1/3}$, (resp. in $V_{2/3}$),
then obviously
$N(\gamma,\mu)\leq N(\VV_1,\mu)$
(resp.
$N(\gamma,\mu)\leq N(\VV_2,\mu)$).
Assume that $\Im\gamma\subset\kr {W}_1$.
 Let $\alpha\in\RRR$
(resp. $\beta\in\RRR$) be the
unique number, such that
$f(\gamma(\alpha)) =2/3-\mu$
(resp.
$f(\gamma(\beta)) =1/3+\mu$).
We have
$\cup_{p\in S(F_i)} U_p(\mu)
\subset
Tb_i(\mu)$
 (for $i=1,2$). Therefore
$\gamma(t)\in\UU (\mu)$
can occur
only if
 $t\geq\beta$ or $t\leq \alpha$.
For $t\leq\alpha $ (resp. $t\geq\beta$)
the curve $\Psi_2^{-1}(\gamma(t))$
(resp. $\Psi_1^{-1}(\gamma(t))$)
 is an integral curve of the vector
field
$(u_2(x),\tau)$
(resp. $(u_1(x),-\tau)$).
Since an integral curve of $u_2$
(resp. $u_1$)
 can intersect no more than
$N(\VV _2,\mu)$
(resp. $N(\VV _1,\mu)$)
standard coordinate neighborhoods of
radius $\mu$, the first part
of (5) follows. The case of curves
in $W_0$ and $W_2$
is considered similarly.

\subsubhead{7. Transversality conditions}
\endsubsubhead

 The next lemma implies that
$\VV$ is almost good
and the (6) of our conclusions.
\proclaim{Lemma 3.4}
\roster \runinitem The family
$\{\big(\cup_{p\in S_i(F_1)} D_p(-u)\big)\}_
{0\leq i\leq n}$
 equals to $\kk (-u_1)$.

{\qquad\qquad}The family
$\{\big(\cup_{q\in S_i(F_2)} D_q(u)\big)\}_
{0\leq i\leq n}$
 equals to $\kk (u_2)$.
\item For each $\lambda\in [0,1]$ the family
$\{\big(\cup_{q\in S_i(F_2)} D_q(-u)\big)
\cap V_\lambda\}_{0\leq i\leq n}$
is an $s$-submanifold of $V_\lambda$,
which is equal to:
\newline
\rm  (1)\it ~
$\Phi (z_1 ,1)(\llll _0(-u_2))\iiif
\lambda =1$,~~
\rm  (2)\it ~
       $\llll _\lambda (-u_2)\iiif
\lambda\in
[1/2, 1 - \epsilon ]$,~~
\rm  (3)\it ~
$\emptyset\iiif
\lambda\leq 1/3$.

\item For each $\lambda\in [0,1]$ the family
$\{\big(\cup_{p\in S_i(F_1)} D_p(u)\big)
\cap V_\lambda\}_{0\leq i\leq n}$
is an $s$-submanifold of $V_\lambda$,
which is equal to:
$$\aligned
&\qquad
(1)~~\Phi (z_0,1)(\llll _0(u_1))\ffor \lambda =0~~;\\
&\qquad
(3)~~\Phi (z_{1/2},1) (\llll _{1/2}(u_1))\ffor
 \lambda =1/2;
\endaligned\qquad
\aligned
&(2)~~\llll _\lambda (u_1)\ffor \lambda\in
[\epsilon,1/2-2\mu ]~~;\\
&(4)~~\emptyset\ffor \lambda\geq 2/3;\\
\endaligned
$$ \item $u$ is an almost good $F$-gradient,
the $s$-submanifold $\kk (u)\cap V_0$
 is almost transversal to $\xx$, and the
s-submanifold  $ \kk (-u)\cap V_1$
 is almost transversal to $\yy$.
\endroster
\endproclaim
\demo{Proof}
 (1) For every $p\in S_i(F_1)$ the positive disc
of $f$ belongs to $V_{1/3}$, which implies
immediately the first
assertion;
the second is proved similarly.

(2)
It is not difficult to see that
it suffices to prove the assertion
for $\lambda\in[-2\mu+2/3,2\mu+2/3]$.
For these values of $\lambda$ it follows from
the analysis of the behaviour of $u$-trajectories
in $\Tb_2(2\mu)$, carried out
in the proof of Subsection 5 above.
 (3) is proved similarly.
(4): To prove that $u$ is almost good let
$p,q\in S(F),~ \ind p\leq \ind q$. Let
$\gamma$ be an $(-u)$-trajectory joining
$p$ with $q$. The case $p,q\in S(F_1)$
or $p,q\in S(F_2)$
follows easily from (1). Let
  $p\in S(F_1), q\in S(F_2)$.
Denote by $z$ the (unique) point
of intersection of
$\gamma$ with
$V_{1/2},~ \ind_Fp$ by $k$,
~ $\ind_Fq$ by $r$. Then
$z\in \big(\Phi(z_{1/2},1)(\llll _{1/2}(u_1))\big)_{k-1}$
and $z\in \llll (-u_2)_{n-r}$,
which is impossible by the choice of $z_{1/2}$.
 The last point is already proved.
\quad$\square$
\enddemo

\document

\subhead{B. Proof of $
\boldkey 1
\boldkey .
\boldkey 1
\boldkey 7
\boldkey (
\boldkey n
\boldkey )
\Rightarrow
\boldkey 1
\boldkey .
\boldkey 1
\boldkey 8
\boldkey (
\boldkey n
\boldkey    +
\boldkey 1
\boldkey )  $
}
\endsubhead
We shall first establish the existence of $\MM$-flows with the
estimate of their quickness
similar to that of 1.18 for
functions on cobordisms without
critical points, then for functions
having one critical point.
The proof will be finished by the usual
induction procedure.
For the rest of \S 3 we fix a natural number
$n$ and we assume that
$1.17(n)$ is true.
We shall need only
the following lemma,
which is an obvious corollary of 1.17(n).

\proclaim{Lemma 3.5} Let $M$ be a closed riemannian
manifold, $\dim M=n$. Let $D>0,~\beta >0,~A>1$.
Then
for every $\mu>0$ sufficiently small
there is an almost good $\MM$-flow $\WW$ such that:
$\text{\rm ~(1)~}$
$\WW$ is of radius $\mu$ and is
$(D,\beta)$-quick;
$\text{\rm ~(2)~}$
 $N(\WW)\leq n2^n$;
$\text{\rm ~(3)~}$
$\GG (\WW )\leq A$. $\square$
          \endproclaim

\subsubhead{1. Functions without critical points}\endsubsubhead

\proclaim{Lemma 3.6} Let $W$ be an $(n+1)$-dimensional
compact cobordism, endowed with a riemannian metric.
Let $g:W\to [a,b]$ be a Morse function without
critical points. Let $C>0$ and let $w$ be an $g$-gradient
such that $\Vert w\Vert\leq C$. Let $\xx$ be an $s$-submanifold
of $V_0$, $\yy$ be an $s$-submanifold of $V_1$.
Let $\mu_0>0$.
Then there is an almost good $\MM$-flow
$\VV =(F,u,\UU)$ on $W$ with the following properties.

 $\text{\rm (1)~}$
  In a
neighborhood of $\partial W$ we have
$u=w$ and $F=g$.
$\text{\rm ~(2)~}$
$D(\VV)\leq\mu_0 ;~~\GG (\VV)\leq 2$.

$\text{\rm  (3)~}$
$\VV$ is $(2C~, 16+n2^{n+5})$-quick, and
$N(\VV)\leq n2^{n+1}$.\qquad
 $\text{\rm ~(4)~}$
$\kk _a(u)\nmid\xx~,~\kk_b(-u)\nmid \yy$.

$\text{\rm  (5)~}$
There are $\alpha~,\beta \in\RRR$, with
$a<\alpha<\beta<b$, such that:
 $g^{-1}(\alpha),~g^{-1}(\beta)$
~et~ $g^{-1}([\alpha,\beta])$ are $(\pm u)$-invariant,
       $g^{-1}([a, \alpha[)$ is $u$-invariant and
$g^{-1}(]\beta, b])$
is $(-u)$-invariant.
For any $u$-trajectory
$\gamma$ in $g^{-1}([a, \alpha[)$
or in $g^{-1}(]\beta, b])$ we have $N(\gamma)\leq n2^n$.

\endproclaim
\demo{Proof} Let $N>0$ be a natural number. For $0\leq s\leq N$
denote $a+{\frac {b-a}N }s$ by $a_s$,
$\frac {2a_i+a_{i+1}}3$ by $b_i$,
 $\frac {a_i+2a_{i+1}}3$
 by $c_i$.
 Denote the cobordism
$g^{-1}([a_s ,a_{s+1}])$ by $W_s$.
Choose $N$ so large that for each $s$
the time which a $\text{\rm grd} (g)$-trajectory
 spends in $W_s$ is not more than
$C/4$. For $0\leq i\leq N-1$ we shall define
by induction in $i$ a sequence of
$\MM$-flows $\VV_i=(F_i,u_i,\UU _i)$ of radius
 $\varkappa_i\leq\mu_0$
on $W_i$.

Let $0\leq i\leq N-1$. For $i>0$ assume that
$\VV_{i-1}=(F_{i-1}, u_{i-1}, \UU_{i-1})$
is already constructed. Denote
$g\vert W_i$ by $g_i:W_i\to[a_i,a_{i+1}]$ and
$w\vert W_i$ by $w_i$. Apply 3.1 to
$g_i, w_i, C$ and get the corresponding
number $\nu_0>0$. Choose
(by 3.5) an $\MM$-flow $\WW_1=(f_1, v_1, \TT_1)$
on $g^{-1}(b_i)$
and
an $\MM$-flow $\WW_2=(f_2, v_2, \TT_2)$
on $g^{-1}(c_i) $,
such that the flow
$\WW_s$ is of radius $\varkappa_i\leq\min(\nu_0,\mu_0)$
and $N(\WW_s)\leq n2^n, \GG(\WW_s)\leq 4/3$,
and $\WW_s$ is $(C/6,1)$-quick
(where $s=1,2$).
Choose the $s$-submanifold of $g^{-1}(a_i)$ as
follows:
$\xx$
for $i=0$ and $\kk(-u_{i-1})\cap g^{-1}(a_i)$
for $i>0$.
Choose the $s$-submanifold of $g^{-1}(a_{i+1})$ as
follows:
empty for $i<N-1$ and $\yy$ for $i=N-1$.
 Theorem 3.1 then
provides
an almost good $\MM$-flow
$\VV_i=(F_i, u_i, \UU_i)$ of radius $\varkappa_i$.
Since $\WW_s$ are $(C/6, 1, \varkappa_i)$-quick,
they are $(C/6, 1+n2^{n+3}, \varkappa_i/2)$-quick
(by 1.13), therefore
$\VV_s$ is
$(2C, 8+n2^{n+4}, \varkappa_i)$-quick.

The property (1) of 3.1 imply that $\VV_i$ glue together
to obtain an $\MM$-flow
$\VV=(F, u, \UU)$ on $W$. I claim
that it satisfies the conclusions of
our Lemma.  (1), (2)
are immediate from (1), (2) of 3.1.
An easy induction argument using
$(\pm u)$-invariance of
$g^{-1}(c_i)$
shows that $\VV$ is almost good.
Proceeding to the estimate of
$N(\gamma)$ and $I(\gamma)$ demanded
by (3) (where $\gamma$ is a $v$-trajectory)
note that the
$(\pm u)$-invariance of
$g^{-1}(b_i)$
and
$g^{-1}(c_i)$
implies easily that $N(\gamma)\leq n2^{n+1}$
and $T(\gamma)\leq 8+n2^{n+4}$
for any $u$-trajectory $\gamma$,
passing by a point of
$g^{-1}([b_i, c_i])$
for some $i$. Assume that $\gamma$
passes by a point of
$g^{-1}(]c_{i-1}, b_i[)$
for some $i$. It is easy to show that
$\gamma$ is defined on $\RRR$.
If $\gamma$ does not intersect
$g^{-1}(a_i)$,
then it stays in $W_i$ or in $W_{i-1}$ and
the required estimates follow (actually
this case does not occur).
If $g(\gamma(t_0))=a_i$, then the
$u_i$-(resp. $(-u_{i-1})$)-invariance
of
$g^{-1}([a_i, b_i[)$
(resp.
$g^{-1}(]c_{i-1}, a_i])$
imply that $g(\gamma(t))\geq a_i$ for $t\geq t_0$
(resp. $t\leq t_0$) and
$N(\gamma)\leq n2^{n+1}, T(\gamma)\leq 2(8+n2^{n+4})$,
and we obtain (3).
Set
$\alpha=b_0, \beta=c_{N-1}$;
then (5) follows from (3) and (5) of 3.1.
Now (4) is easy to prove.
~$\square$
\enddemo

\subsubhead{2. Functions with one critical point}
\endsubsubhead

\proclaim{Lemma 3.7}
Let $W$ be an $(n+1)$-dimensional compact
riemannian
cobordism.
Let $g:W\to [a,b]$ be a Morse function
with one critical point $p$ and a strongly standard chart
$\Phi_p:U_p\to B^{n+1} (0,r_p)$
 around $p$.
Let $w$ be a $g$-gradient with respect to
$\{\Phi_p\}$.
Then for every $\delta>0$ sufficiently small
 the time which a $w$-trajectory
can spend in
$g^{-1}([g(p)-\delta^2,~g(p)+\delta^2])
\setminus U_p(\delta/2)$ is not more than
$6$.
\endproclaim
\demo{Proof}
 Denote
$ g^{-1} ([g(p)-\delta^2,g(p)+\delta^2])$ by $S_\delta $.
Let $\gamma$ be a $w$-trajectory.
The function $dg(w)$ is bounded from below in
$S_\delta
\setminus U_p$, therefore for
$\delta >0$ sufficiently small the time which
$\gamma$ spends in $S_\delta\setminus U_p$
 is less
than 1. The time, which $\gamma$ spends in
$U_p\cap \big(S_\delta
\setminus U_p(\delta)\big)$
is not more than 2 (by 1.12), and
 the time which $\gamma$ spends
 in $U_p(\delta)\setminus U_p(\delta/2)$
is not more than
$\ln 8\leq 3$ (by 1.11).
\quad$\square$
\enddemo

\proclaim{Lemma 3.8}
Let $W$ be an $(n+1)$-dimensional compact
cobordism endowed with a riemannian metric.
Let $(g,w,\UU)$ be an $\MM$-flow on $W$, where
 $g:W\to [a,b]$ is a Morse function
with one critical point $p,~ g(p)=c$,
and the $g$-chart
around $p$ is strongly standard.
Let $C>0$, and assume that $\Vert w\Vert\leq C$.
Let $\mu_0>0$.
Let $\xx$ be an $s$-submanifold of $V_0$, and
$\yy$ be an $s$-submanifold of $V_1$.
Then there is an almost good $\MM$-flow
$\VV=(F,u,\UU)$ on $W$,
with the following properties.

 $\text{\rm (1)~}$
 In a
neighborhood of $\partial W$ we have
$u=w$ and $F=g$.
$\text{\rm ~(2)~}$
$D(\VV)\leq\mu_0; \GG (\VV)\leq 2$.

$\text{\rm  (3)~}$
$\VV$ is $(2C~, 40+n2^{n+6})$-quick, and
 $N(\VV)\leq 1+n2^{n+1}$.\quad
 $\text{\rm ~(4)~}$
$\kk _a(u)\nmid\xx~,~\kk_b(-u)\nmid \yy$.

$\text{\rm  (5)~}$
There are $a'~,b' \in\RRR$ with
$a<a'<c<b'<b$ such that:
 $g^{-1}(a'),~g^{-1}(b'),
\mathbreak
~g^{-1}([a',b'])$
are $(\pm u)$-invariant;
 $g^{-1}([a, a'[)$ is $u$-invariant,
$g^{-1}(]b', b])$
is $(-u)$-invariant.
 For any $u$-trajectory
$\gamma$ in $g^{-1}([a, a'[)$
or in
$g^{-1}(]b', b])$ we have $N(\gamma)\leq n2^n$.
\endproclaim

\demo{Proof}
Let $\delta\in ]0,\mu_0[$ be so small that
 it is less than the radius of the standard
chart $\Phi_p$ of $g$ around $p$, and that the
conclusions of 3.7 hold and that
$\GG (U_p(\delta),\Phi_p)\leq 2$.
Denote $g^{-1}(c-\delta^2)$ by $V_- ~$,
$g^{-1}([a, c-\delta^2])$ by $W_- ~$,
$g^{-1}(c+\delta^2)$ by $V_+ ~$,
$g^{-1}([c+\delta^2, b])$ by $W_+ ~$,
$g^{-1}([c-\delta^2, c+\delta^2])$ by $W'$.
Apply Lemma 3.6 to the function $g\vert W_+$,
$(g\vert W_+)$-gradient $w$, the number $\mu_0$, the
$s$-submanifold $\yy$ of
$V_1$,  and the $s$-submanifold
$D(p,-w)\cap V_+$ of $V_+$
(consisting of one compact manifold). Denote
the resulting flow by $\VV_+=(F_+, u_+, \UU_+)$.
The corresponding values $\alpha,\beta$,
satisfying (5) of 3.6 with respect to
$\VV_+$ will be denoted by $b'',b'$; we have $b''<b'$.
Denote by $\Phi'_p:U'_p\to B^{n+1}(0,\delta/2)$ the restriction
of the $g$-chart around $p$, and
by $\VV'$ the $\MM$-flow $(g\vert W',w\vert W',
\Phi'_p)$.
 (1) of 3.6 implies that we can glue
$\VV_+$ to $\VV'$. Denote the resulting $\MM$-flow by
$\VV_1=(F_1,~u_1,~\UU_1)$. It is almost good since
$\kk(u_+)\cap
V_+\nmid D(p,-w)\cap V_+$.
Apply Lemma 3.6 to $g\vert W_-$, the
$(g\vert W_-)$-gradient $w$, the number $\mu_0$,
the $s$-submanifold $\kk (u_1)\cap V_-$ of $V_-$ and to the
$s$-submanifold $\xx$ of $V_0$.
Denote the resulting $\MM$-flow by
$\VV_-=(F_-,u_-,\UU_-)$. The corresponding values
$\alpha,~\beta$ satisfying (5) of 3.6
with respect to $\VV_-$ will be denoted by $a',~a''$, so that $a'<a''$.
 (1) of 3.6 implies that we can glue
$\VV_-$ to $\VV_1$. Denote the resulting $\MM$-flow
on $W$ by $\VV=(F,u,\UU)$. It is almost good, since
$\kk(u_1)\cap
V_-\nmid \kk (-u_-)\cap V_-$.

We claim that the $\MM$-flow $\VV$ together with the
numbers $a',b'$ satisfy the conclusions of our Lemma.
Indeed, (1),(2) and (5) follow from the construction.
(4) follows since
$g^{-1}(]b',b])$ (resp. $g^{-1}([a,a'[)$)
is $u$-invariant and from the
corresponding transversality property
of $\VV_+$  (resp. of $\VV_-$).
To prove (3) note that $(\pm u)$-invariance
of $g^{-1}(a'')$ and of
$g^{-1}(b'')$
together with (3) of 3.6  imply
already the required estimates of $N(\gamma)$ and of $T(\gamma)$
for all the $u$-trajectories
$\gamma$
 passing by a point of
$g^{-1}([a,~a''])$
 or of $g^{-1}([b'',b])$.
Let now $\gamma$ be an $u$-trajectory, passing by a point of
$g^{-1}(]a'',~b''[)$.
Since $g^{-1}(a'')$ and $g^{-1}(b'')$
are $(\pm u)$-invariant, $\gamma$ is defined
on $\RRR$. If $\gamma$ does not intersect
$W'$, then
$\Im\gamma\subset W_+$
or
$\Im\gamma\subset W_-$
and the required
estimates of $T(\gamma)$ and $N(\gamma)$
again follow from those of 3.6 (actually
this case does not occur). If there is $t_0\in\RRR$
such that $\gamma(t_0)\in g^{-1}([c-\delta^2, c+\delta^2])$,
then it is not difficult to show
(using the fact that $u$ is a $g$-gradient
in $W'$,  $u$-invariance of
$g^{-1}([c+\delta^2, b''[)$ and
$(-u)$-invariance of
$g^{-1}([a'', c-\delta^2]$)
that three possibilities can occur for $\gamma$:

 $\text{\rm(1)~}$
 $\exists\alpha,\beta\in\RRR:
~\gamma(]-\infty,\alpha])\subset W_-,\quad
\gamma([\beta, \infty[)\subset W_+,\quad
\gamma([\alpha,\beta])\subset W'$;\newline
 $\text{\rm~(2)~}$
      $\exists\beta\in\RRR:
~\gamma([\beta, \infty[)\subset W_+ ,\quad
\gamma(]-\infty,\beta])\subset g^{-1}(]c, c+\delta^2])$, and
$\lim_{t\to-\infty}\gamma(t)=p$; \newline
 $\text{\rm~(3)~}$
$\exists\alpha\in\RRR:~ \gamma(]-\infty,\alpha])\subset W_-,\quad
\gamma([a,\infty[)\subset g^{-1}([ c-\delta^2,c[)$,
and
$\lim_{t\to\infty}\gamma(t)=p$.

It is easy to check the required estimates
for $T(\gamma)$ and $N(\gamma)$ in each of these cases.
~$\square$\enddemo

\subsubhead{3. End of the proof of
 1.17(n)~$
\Rightarrow$
1.18(n+1)}\endsubsubhead

We are given $B>0,\mu_0>0$.
Let $\VV_0=(g,w,\UU)$ be an $\MM$-flow
on $M$,
such that
 $\Vert w\Vert\leq B/2$,
 all the charts of $\UU$ are strongly standard, and
 if $p,q\in S(g)$
and $p\not= q$, then $g(p)\not= g(q)$
(existence of such a flow
is an easy exercise in Morse theory).
Let $p_1,...,p_r$ be the critical points of $g$, so that
$g(p_i)<g(p_{i+1})$. Denote $g(p_i)$ by $a_i$, and set
$b_i=\frac {a_i+a_{i+1}}2$ for $1\leq i\leq r-1$,
$b_0=a_1-1,~b_r = a_r+1$. Denote $g^{-1}([b_i ,b_{i+1}])$
by $W_i$. We shall construct by induction in $i$
a sequence of $\MM$-flows $\VV_i=(f_i,v_i,\UU_i)$
on $W_i$.

Let $0\leq i\leq r-1$. For $i>0$ assume that
the $\MM$-flow $\VV_{i-1}$ is already constructed.
Denote $g\vert W_i$ by $g_i:W_i\to[b_i,b_{i+1}]$,
$w\vert W_i$ by $w_i$. Choose some restriction
$\UU'_i$ of the standard chart for $g_i$
to obtain an $\MM$-flow
$\WW_i=(g_i, w_i, \UU'_i)$ on $W_i$;
we have $\Vert w_i\Vert\leq B/2$.
For $i=0$ choose the
empty $s$-submanifold
of $f^{-1}(b_0)=\emptyset$
and the empty $s$-submanifold of $f^{-1}(b_1)$.
For $0<i$ choose the $s$-submanifold
$\kk(-u_{i-1})\cap f^{-1}(b_i)$
of  $f^{-1}(b_i)$
and the empty $s$-submanifold of
$f^{-1}(b_{i+1})$.
Applying
3.8, we obtain an $\MM$-flow $\VV_i$ on $W_i$.
The corresponding numbers $a', b'$ will be denoted by
$c_i, d_i$; then $b_i<c_i<a_{i+1}<d_i<b_{i+1}$.
 (1) of 3.8 implies that $\VV_i$
glue together to an $\MM$-flow on $M$, which will be denoted
by $\VV=(f,v,\UU)$. An easy induction argument
using the $(\pm v)$-invariance of
$g^{-1}(c_i)$ and the almost transversality
of $\kk (-v_i)\cap g^{-1}(b_{i+1})$ to
$\kk (v_{i+1})\cap g^{-1}(b_{i+1})$
shows that $\VV$ is almost good. The property
(1) of 1.18 follow immediately from (2) of 3.8.
The estimates of $T(\gamma)$ and $N(\gamma)$
demanded by (2) and (3) of 1.18 follow from
the $(\pm v)$-invariance of $g^{-1}(c_i), g^{-1}(d_i)$
similarly to the estimates of $T(\gamma)$ and
$N(\gamma)$ of Lemma 3.6.
\quad$\square$

\head{\S4. Ranging systems and proof of  Main Theorem} \endhead

In \S4          the excision
isomorphisms
$H_*(X\setminus B,
 A\setminus B)\to H_*(X, A)$  are denoted
 by Exc.
\subhead{A. Generalities on intersection indices}\endsubhead
Let $M$ be a manifold without boundary,
 $\dim M =m$.
Let $X\subset M$ and $N$ be an oriented submanifold
of
$M$, such that $N\setminus \Int X$ is compact.
Then the orientation class
$\mu_{N\setminus\Int X}\in
H_n(N,N\cap\Int X)$
is defined, where $n=\dim N$
(see [7, Th.A8]). Its image in $H_n(M,X)$ will be
denoted by $[N]_{M,X}$ (or simply by $[N]$ if there is no
possibility of confusion).
The next lemma follows immediately from the theorem cited above.

\proclaim{Lemma 4.1}
\roster\runinitem Let $Y$ be a closed subset of $\Int X$. Then
$[N\setminus Y]_{M\setminus Y,X\setminus Y}$
equals to
${\qquad\quad}$
the image of
$[N]_{M,X}$ with respect
to $\text{\rm Exc }^{-1}
 : H_*(M,X)\overset{\approx}
\to\longrightarrow
H_*(M\setminus Y,X\setminus Y)$.

\item Let $M'$ be a manifold, $U$ be an open subset
of $M'~,~X'\subset U$.
Let $\phi:M\to
U$ be a diffeomorphism such that
$\phi (X)\subset X'$. Denote $\phi (N)$ by $N'$.
Denote by $\Phi:(M,X)\to (M',X')$
the resulting map of pairs.
Then $[N']_{M',X'}=\Phi_
 *([N]_{M,X}).~\square$
\endroster
\endproclaim

 Let $L$ be a compact cooriented submanifold
without boundary
of
 $M$.
Then the canonical coorientation class
 $]L[\in H^{m-l}(M, M\setminus L)$
is defined, where $l=\dim L$.
Assume that $X\cap L=\emptyset ~,~ N\pitchfork L \aand n+l=m$.
Denote by $j$ the embedding $(M,X)\hookrightarrow (M,M\setminus L)$.
The image of the class $[N]$ in $H_n(M, M\setminus L)$ will be
denoted again
by $[N]$. The following lemma is standard.

\proclaim{Lemma 4.2}
In the above assumptions the set $N\cap L$ is finite and
the intersection index $N ~\sharp ~ L$ equals
 $\big\langle j^*(]L[) , [N]\big\rangle. \quad\square$
\endproclaim

\subhead{B. Ranging systems and $\st v$-images of the
fundamental classes}\endsubhead
Let $f:W\to [a,b]$ be a Morse function on a compact
riemannian
cobordism $W$,
$v$ be an $f$-gradient, and denote
 $\fpr fa$ by $V_0$,$\fpr fb$ by $V_1$.
This terminology is valid for Subsections B and C.

\definition{Definition 4.3}
Let $\Lambda =\{\lambda_0,...,\lambda_k\}$ be a finite set of
regular values of $f$, such that $\lambda_0 =a,\lambda_k =b$,
and for each $0\leqslant i\leqslant k-1$ we have $\lambda_i <\lambda_{i+1}$
and there is exactly one critical value of $f$ in
$[\lambda_i ,\lambda_{i+1}]$. The values $\lambda_i ,\lambda_{i+1}$
will be called adjacent. The set of pairs $\ran$ is called
\it ranging
system
for $(f,v)$
\rm
 if
\roster
\item"(RS1)" $\forall \lambda\in \Lambda$ $A_\lambda$ and
$B_\lambda$ are disjoint compacts in $\fpr f{\lambda}$.

\item"(RS2)" Let $\lambda,\mu\in\Lambda$ be adjacent. Then for
 every $p\in S(f)\cap \fpr f{[\lambda,\mu]}$ either

i) $D(p,v)\cap \fpr f{\lambda}\subset \Int A_\lambda$
\qquad
 or \qquad
ii) $D(p,-v)\cap \fpr f{\mu}\subset \Int B_\mu$.

\item"(RS3)"  Let $\lambda,\mu\in \Lambda$ be adjacent. Then
$\stind v\mu\lambda (A_\mu)\subset \Int A_\lambda ~\text{and}~
\stind {(-v)}\lambda\mu (B_\lambda)\subset
\Int B_\mu.
 \triangle$
\endroster
\enddefinition

 The following
 lemma  is obvious.
\proclaim{Lemma 4.4}
\roster
\runinitem Let $S=\ran$ be a ranging system for $(f,v)$.
Let $\mu,\nu\in\Lambda$.
${\qquad\quad}$
 Then
$S'= \{(A_\lambda , B_\lambda)\}_{\lambda\in\Lambda,
\mu\leqslant\lambda\leqslant\nu}$ is a ranging system for
$\big(f\vert \fpr f{[\mu,\nu]},v)$.
\item Let $c\in ]a,b[$ be a regular value of $f$.
Let $\{(A'_{\lambda'} , B'_{\lambda'})\}_{\lambda'\in\Lambda'}$,
\break
$\{(A''_{\lambda''} , B''_{\lambda''})\}_{\lambda''\in\Lambda''}$
be ranging systems for $(f\vert \fpr f{[a,c]} ,v)$, and,
 respectively, for
$(f\vert \fpr f{[c,b]} ,v)$. Assume that
$A''_c\subset A'_c ~,~   B'_c\subset B''_c  \aand A'_c\cap
B''_c=\emptyset$.
Then the system $\ran$,
defined by (C) is a ranging system for $(f,v)$.
$$
\left.
\gathered
\Lambda =\Lambda'\cup\Lambda'';\qquad\qquad
A_\lambda =A'_{\lambda}\aand B_\lambda =B'_{\lambda}
\ffor\lambda\in\Lambda'\setminus\{c\} \\
A_\lambda =A''_\lambda\aand B_\lambda =B''_\lambda
\ffor\lambda\in\Lambda''\setminus\{c\};\quad
A_c =A'_c,\quad
B_c=B''_c
\endgathered \quad \right\}
\text{(C)}
\gathered
{}\\
\quad \square
\endgathered
$$
\endroster
\endproclaim
In 4.5-4.7 $\ran$ stands for a ranging system for $(f,v)$.

\proclaim{Lemma 4.5}
There is $\epsilon>0$ such that for any $f$-gradient $w$ with
$\nr v-w\nr <\epsilon$ the ranging system $\ran$ is also
a ranging system for $(f,w)$.
\endproclaim

\proclaim{Proposition 4.6}
Let $N$ be a submanifold of $V_1\setminus B_b$
   such that
$    N\setminus \text{\rm Int}~
    A_b$ is compact.

Then $N'=\stind vba (N)$ is a submanifold of $V_0\setminus B_a$
  such that
$     N' \setminus \text{\rm Int}~ A_a$     is compact.
\endproclaim

\proclaim{Proposition 4.7}
There is a homomorphism
$\qquad H(v):H_*(V_1\setminus B_b,A_b)\longrightarrow
\break
H_*(V_0\setminus B_a,A_a)$, such that
\roster
\item If $N$ is an oriented submanifold of $V_1$, satisfying the
hypotheses of 4.6, then
\newline
$H(v)([N])=[\stind vba (N)]$.

\item There is an $\epsilon>0$ such that for
 every $f$-gradient $w$ with $\nr w-v\nr <\epsilon$
we have $H(v)=H(w)$.

\endroster
\endproclaim

The proof of 4.5 - 4.7 occupies the rest of
Subsection B.
An easy induction argument shows that it is
 sufficient to
prove each of them in the case $\text{card}~\Lambda =1$.
Let $S(f)=S1(f)\sqcup S2(f)$, where
for every $p\in S1(f)$, resp. $p\in S2(f)$
the i), resp. ii) of (RS2) holds.
Pick Morse functions $\phi_1, \phi_2:W\to[a,b]$, adjusted
to $(f,v)$, such that there are regular values
$\mu_1$ of $\phi_1$, $\mu_2$ of $\phi_2$
satisfying:
(1) for every $p\in S1(f)$ we have: $\phi_1(p)<\mu_1$
and $\phi_2(p)>\mu_2$.
(2) for every $p\in S2(f)$ we have: $\phi_1(p)>\mu_1$
and $\phi_2(p)<\mu_2$.
\subsubhead
Proof of 4.5
\endsubsubhead
Proposition 1.10 implies that (RS2) holds for
all $f$-gradients $w$, sufficiently close to $v$.
Passing to (RS3), consider
the cobordism $\phi_1^{-1}([\mu_1,b])$ and the
$\phi_1$-gradient $v$. Denote
$T(A_b,v)\cap\phi_1^{-1}(\mu_1)$
by $Z_+$.
Since every $(-v)$-trajectory starting in
$A_b$ reaches $\phi_1^{-1}(\mu_1)$, ~$Z_+$ is compact.
Consider the cobordism $\phi_1^{-1}([a,\mu_1])$
and the compact        $V_0\setminus\Int A_a\subset V_0$.
Every $v$-trajectory starting in
$V_0\setminus\Int A_a$ reaches $\phi_1^{-1}(\mu_1)$
and  the compacts $Z_+$ and
$Z_-=\phi_1^{-1}(\mu_1)\cap T(V_0\setminus \Int A_a,-v)$
are
 disjoint.
Choose disjoint open neighborhoods:
$U_+$ of $Z_+$ and $U_-$ of
$Z_-$ in $\phi_1^{-1}(\mu_1)$.
It follows from 5.6 that
there is $\epsilon''>0$
such that for every $w\in\vew$
 with $\Vert w-v\Vert<\epsilon''$
we have:
$$
T(A_b ,w)\cap\phi_1^{-1}(\mu_1)\subset U_+,\aand
T(V_0\setminus\Int A_a,-w)\cap\phi_1^{-1}(\mu_1)
\subset U_-,
$$
therefore for all $f$-gradients $w$,
sufficiently $C^0$ close to
$v$, we have:
$\st w (A_b)\subset \Int A_a$. The second part of
(RS3) is considered in the same way.~$\square$

For $\delta>0$ denote by $D1_\delta (v)$, resp. by
$D1_\delta(-v)$,  the intersection
with $V_0$, resp. with $V_1$,
of $\cup_{p\in S1(f)} D_\delta (p,v)$,
resp. of $\cup_{p\in S1(f)} D_\delta (p,-v)$.
By abuse of notation
the intersection of
$\cup_{p\in S1(f)} D(p,v)$
with $V_0$ will be denoted by
$D1_0(v)$.
Denote $D2_\delta(-v)\cup\stexp {(-v)} (B_a)$ by $\Delta(\delta,-v)$
and
$D1_\delta(v)\cup\st v (A_b)$ by $\nabla(\delta,v)$.
The similar notation like
$D2_\delta(-v),~ \Delta(0,-v),
\mathbreak
    D2_0(v)  $
etc.
are now clear without special definition.
For
 $\delta>0$    sufficiently small  we have
$$\gathered
\forall p\in S1(f): D_\delta(p)\in\phi_1^{-1}(]a,\mu_1[)\aand D_\delta(p)\subset
\phi_2^{-1}(]\mu_2, b[)\\
\forall p\in S2(f): D_\delta(p)\in\phi_1^{-1}(]\mu_2, b[)\aand D_\delta(p)\subset
\phi_2^{-1}(]a, \mu_2[)
\endgathered\tag\text{D1}
$$
Applying 2.9 and 2.10
to functions $\phi_1,\phi_2$ and their restrictions
 it is easy to prove that for $\delta>0$
sufficiently small we have:
$$\nabla(\delta,v)\subset \Int A_a,\quad
\Delta(\delta,-v)\subset\Int B_b,\quad
\Delta(\delta,-v)\cap
D1_\delta(-v)=\emptyset
\tag\text{D2}
$$
Fix some $\delta>0$ satisfying (D1) and (D2).

\subsubhead
{Proof of 4.6}
\endsubsubhead
 The set $N'\setminus \Int A_a$ is a closed
subset of
$\st v \big( N\setminus
     (\Int A_b\cup B_\delta (-v)    )\big) $.
The set $ N\setminus
\big(\Int A_b\cup B_\delta (-v)\big)$
is a compact subset of the domain of $\st v$
and 4.6 follows. $\qquad\square$

\subsubhead{  Homomorphism
$H(v;\mu',\mu; U): H_*(V_1\setminus B_b,A_b)\to
H_*(V_0\setminus B_a,A_a)$}
\endsubsubhead

Let $0\leqslant\mu '<\mu\leqslant\delta$.
Let $U$ be any subset of $V_1$ such that
$$\Delta (0,-v)\subset U
\subset B_b
\aand U\cap D1_{\delta}(-v)=\emptyset
\tag4.1$$
(for example $U=\Delta (\delta ,-v)$ will do).
 Denote by
$\HH (v;\mu',\mu ; U)$
the following sequence of homomorphisms
$$\gather
H_*(V_1\setminus B_b,A_b)
@>        I_*   >>
H_*\big( V_1\setminus U~,~
A_b\cup D1_\mu (-v))
@> {\text{Exc} ^{-1}} >>\\
H_*\big( V_1\setminus (U
\cup D1_{\mu'}(-v)\big)~,~
\big(A_b\cup D1_\mu (-v)\big)\setminus D1_{\mu'}(-v)\big)
@> \st v_* >>
  H_*(V_0\setminus B_a,A_a).
\endgather
$$
(Here $I$
is the corresponding inclusion.  Note that the
 last arrow is well defined
since
$\stexp {(-v)} (B_a)\subset U\aand D_0(-v)\cap V_1\subset
D2_0(-v)\cup D1_{\mu'}(-v).$)
 The composition
$\st v_* \circ {\text{Exc} ^{-1}} \circ I_* $
of this sequence
will be denoted by $H(v; \mu',\mu;U)$.
\subsubhead{  Homomorphism $
H(v): H_*(V_1\setminus B_b,A_b)\to
H_*(V_0\setminus B_a,A_a)$}
\endsubsubhead

Suppose that  $U'$
is another subset of $B_b$
 satisfying $(4.1)$ and $U\subset U'$. Then
the inclusion $V\setminus U'\subset V_1\setminus U$
induces a map of the sequence
$\HH (v;\mu',\mu ; U')$ to the sequence
$\HH (v;\mu',\mu ; U)$,
which is identity on the first term
$H_*(V_1\setminus B_b,A_b)$
and on the last term
$H_*(V_0\setminus B_a,A_a).$
This implies easily that
$   H(v;\mu',\mu ; U)$ does not depend on the choice
of $U$. Similarly one checks that
$   H(v;\mu',\mu ; U)  $ does not depend
on the choice of $\mu'$ and $\mu$,
neither on the choice of
$\delta$, or on the choice
of presentation $S(f)=S1(f)\cup S2(f)$
(if there is more then one such presentation).
Therefore this homomorphism is well determined by
$v$ and the ranging system $\ran$.
We shall denote it by $H (v)$.

\subsubhead{Proof of (1) of 4.7}
\endsubsubhead
Consider the sequence
$\HH (v;0,\delta; U)$ with any $U$ compact
and apply Lemma 4.1.\quad$\square$

Passing to the proof of  4.7 (2)
choose and fix some  $\mu',\mu$ with
$0<\mu'<\mu<\delta$ and an open subset $U$,
such that
$
\Delta(0,-v)\subset U\subset\Int B_b,\aand D1_\delta(-v)
\subset V_1\setminus\overline U.
$
The proof of the next lemma is similar to the proof
of 4.5 above and will be omitted.

\proclaim{Lemma 4.8}
There is $\epsilon>0$ such that for every
$f$-gradient $w$ with $\Vert w-v\Vert<\epsilon$
we have:
$$
\gather
\Delta(0,-w)\subset U;\quad \nabla(\delta,w)\subset\Int A_a\tag4.2\\
D1_0(-w)\subset B1_{\mu'}(-v);
\quad D1_\mu(-v)\subset B1_\delta(-w);\quad D1_\delta(-w)
\subset V_1\setminus\overline U.\qquad\square\tag4.3
\endgather
$$
\endproclaim
\demo{Proof of (2) of 4.7}
Denote $V_1\setminus \big(U\cup B1_{\mu'} (-v)\big)$
by $R$ and
$\big(A_b\cup D1_\mu (-v)\big)\setminus B1_{\mu'}(-v)\big)$
by $Q$.
Then $R$ is in the domain of
 $\st v$, and $\st v$ maps
the pair of compacts
$(R, Q)$
to the pair of open sets $(V_0\setminus B_a~,~ \Int A_a)$.
By 5.6 there is $\kappa>0$
such that for every
$u\in\vew$ with $\Vert u-v\Vert<\kappa$
each $(-u)$-trajectory
starting at a point $x\in R$
reaches $V_0$ and intersects it at a point
$L(u,x)\in V_0\setminus B_a$, and  for
$x\in Q$
we have $L(u,x)\in\Int A_a$.
By [19, 8.10] the map $(u,x)\mapsto L(u,x)$
is continuous with respect to  $C^1$ topology in
$\vew$.
In particular for every $f$-gradient $w$ with $\Vert w-v\Vert<\kappa$
\quad $\st w$ defines a map of pairs
\break
$(R,Q)
\to
(V_0\setminus B_a~,~ \Int A_a)$.
This map is homotopic to $\st v$ by the homotopy
${\Cal H}_t(x)=L(tv+(1-t)w, x )$. Therefore $H(v)$
equals to
the composition of the following
sequence:
$$
\gathered
H_*(V_1\setminus B_b,A_b)
@> I_*          >>
H_*\big( V_1\setminus U~,~
A_b\cup D1_\mu (-v)\big)
@> {\text{Exc} ^{-1}} >> \\
H_*\big( V_1\setminus (U
\cup D1_{\mu'}(-v)\big)~,~
\big(A_b\cup D1_\mu (-v)\big)\setminus D1_{\mu'}(-v)\big)
@> \st{ w}_* >>
  H_*(V_0\setminus B_a,A_a).
\endgathered
\tag4.4
$$
Let $\epsilon>0$ satisfy the conclusions
of 4.8, and set $\epsilon_0=\min (\epsilon,\kappa)$.
Let $w$ be an $f$-gradient
with $\Vert w-v\Vert<\epsilon_0$.
Then (4.2), (4.3) imply that
there is a morphism of
 the sequence
(4.4) to the
sequence $\HH (w; 0,\delta; U)$ such
that the homomorphisms on the first term
and on the last term
are identical (the two homomorphisms in the middle
are induced by inclusions). This implies that the
composition of (4.4) equals to $H(w)$ and 4.7 follows.\quad
$\square$
\enddemo

\subhead{C. Constructing ranging systems}\endsubhead
In this subsection we construct ranging
systems from $\MM$-flows
on $V_0$ and $V_1$.
 Denote $\dim W$ by $n$. We assume that $S(f)\not=\emptyset$.
\definition{Definition 4.9 }
Let $v$ be an almost good $f$-gradient. Let
$\VV _0=
\break
(\phi_0,u_0,\UU _0),
\VV _1=(\phi_1,u_1,\UU _1)$ be almost good $\MM$-flows on $V_0$,
resp.$V_1$.
We shall say that $(\VV _0,\VV _1)$ is a
\it ranging pair
for $(f,v)$ \rm
if
$\kk_a(v)\nmid\kk(-u_0),~
\kk_b(-v)\nmid\kk(u_1)$, and
there is a $\delta>0$, such that for any
 $0\leqslant s\leqslant n-1$ the
following conditions hold:
$$
\gather
\text{ for some }~~
\delta_1>\delta~~
\text{the gradients}~~
v,u_0\aand u_1 \quad \text{are}\quad \delta_1 -
\text{separated}
\tag\text{RP}1\\
\overline{\TTT (\KKK (u_1), v) _{(\leqslant s+1)}
 (\delta)} \cap V_0
\subset
\KKK (u_0)_
{(\leqslant s)}(\delta)\tag\text{RP}2
\\
\overline{\TTT    (\KKK (-u_0), -v )
_{(\leqslant s+1)}     (\delta)}
\cap V_1
\subset
\KKK (-u_1)_{(\leqslant s)}(\delta)
\qquad\qquad\triangle \tag\text{RP}3
\endgather
$$
\enddefinition
\remark{Remark }
Conditions (RP2) and    (RP3)
demand  in particular
that $\delta$ is in
the interval of definition of
 $\TTT(\KKK(u_1),v)$,
resp.
$\TTT(\KKK(-u_0),-v)$
\quad$\triangle$
\endremark

We shall show that every
 ranging pair $(\VV _0,\VV _1)$ generates
a bunch
 of ranging
 systems.

\definition{Construction 4.10}
Let $\Lambda=\{\lambda_0,...\lambda_k\}$ be a set of
 regular values of $f$, such
that $\lambda_i<\lambda_{i+1},\lambda_0=a,\lambda_k=b$
 and in each $[\lambda_i,\lambda_{i+1}]$
there is only one critical value of $f$.
 Choose some $\delta'\in ]\delta, \delta_1[$
 in such a way, that
$\forall s: 0\leqslant s\leqslant n-1$ we have
$$
\gathered
\overline{\TTT (\KKK (u_1), v) _{(\leqslant s+1)}
 (\delta')} \cap V_0
\subset
\KKK (u_0)_
{(\leqslant s)}(\delta)\\
\overline{\TTT    (\KKK (-u_0), -v )
_{(\leqslant s+1)}     (\delta')}
\cap V_1
\subset
\KKK (-u_1)_{(\leqslant s)}(\delta)
\endgathered
\tag\text{RP}4
$$
(such $\delta'$
 exists since $\delta$ is in the interval of definition
of $\TTT(\KKK(u_1),v)$, resp.    $\TTT(\KKK(-u_0),-v)$.)

We shall now define for each integer
$s:0\leqslant s\leqslant n$ compact subsets
$A_\lambda^{(s)}$ and $B_\lambda^{(s)}$
 of $\fpr f\lambda$.
Set $A_\lambda^{(0)} = B_\lambda^{(0)} =\emptyset$.
Let $\delta=\delta_k<\delta_{k-1}<...<\delta_1<\delta_0=\delta'$
 be a sequence of real numbers. For  $1\leqslant s\leqslant n$ set
$$
\align &
\left\{
\aligned     A_{\lambda_l}^{(s)}  &=
\overline{\TTT (\KKK (u_1), v) _{(\leqslant s)}
     (\delta_l)} \cap f^{-1}(\lambda_l)
\ffor 0< l\leqslant k \\
 A_a^{(s)}  &= ~~ \overline{\KKK (u_0)_{(\leqslant s-1)}(\delta)}
 ~~= D_{\delta} (\indl {s-1};u_0) \endaligned
\right.\\
&
\left\{
\aligned     B_{\lambda_l}^{(s)}  &=
\overline{\TTT    (\KKK (-u_0), -v )
_{(\leqslant s)}     (\delta_{k-l})}\cap f^{-1}(\lambda_l)
\ffor 0\leqslant l< k \\
 B_b^{(s)}  &= ~~ \overline{\KKK (-u_1)_{(\leqslant s-1)}(\delta)}
 ~~= D_{\delta} (\indl {s-1};-u_1) \endaligned
\right.
\endalign
$$
Lemma 2.16 imply that
$A_b^{(s)} = D_{\delta} (\indl {s-1};u_1)\aand B_a^{(s)} =
D_{\delta}(\indl {s-1};-u_0) $. \quad$\triangle$
\enddefinition
\proclaim{Lemma 4.11}
\roster
\runinitem $A_\lambda^{(r)}\cap B_\lambda^{(s)}=\emptyset$
for $r+s\leqslant n$.
\item If $0\leqslant s\leqslant n-1$, then
$\{(A_\lambda^{(s)}, B_\lambda^{(n-s-1)})\}_{\lambda\in\Lambda}$ is
 a ranging system for $(f,v)$.
\item For every $0\leqslant s\leqslant n-1$ the pair
$(V_1\setminus B_b^{(n-s-1)} ,A_b^{(s)})$ is homotopy
equivalent to a finite CW-pair,
 having only cells of dimension $s$.
\endroster
\endproclaim
\demo{Proof}
(1) For $\lambda =a,\lambda =b$ it follows directly
from (RP1). Let
$0<l<k$ and set $\lambda =\lambda_l$.
Then $\delta_l~,~\delta_{k-l} <\delta'\aand
 A_\lambda ^{(r)}
\subset
T(\KKK (u_1)_{(\leqslant r-1)}(\delta'),v)
\cup B_{\delta'}(\indl r;v),\aand
B_\lambda^{(s)}\subset
T(\KKK (-u_0)_{(\leq s-1)}(\delta'),-v)
\cup B_{\delta'}(\indl s ,-v)$.
Therefore if $A_\lambda ^{(r)}\cap B_{\lambda}^{(s)}\not=\emptyset$,
there is a $(-v)$-trajectory
joining
\quad
1) a point of $\KKK (u_1)_{(\leqslant r-1)}(\delta')$ with a point
 of $\KKK (-u_0)_{(\leqslant s-1)}(\delta')$ or

2) a point of $\KKK (u_1)_{(\leqslant r-1)}(\delta')$
with a point of
$B_{\delta'} (p)$,~where~  $p\in S(f),\ind p\geq n-s$
or

3) a point of $\KKK (-u_0)_{(\leqslant s-1)}(\delta')$
 with a point of
$B_{\delta'} (q),\wwhere q\in S(f),\ind q\leqslant r$
or

4) a point of
$B_{\delta'} (p)$ with a point of
 $B_{\delta'} (q)$, where $p,q\in S(f),\ind q\geq n-s,\ind p\leq r$.

The last case        is impossible,
 since $v$ is
$\delta'$-separated, the first is impossible
because of (RP4) and since $u_0$ is
$\delta'$-separated. 2) and 3) are considered similarly.

(2) (RS1) follows from (1). (RS2) follows, since for all
$p\in S(f)$ we have $\ind p\leqslant s$ or $\ind p\geq s+1$.
 (RS3) follows from the construction.

(3) Let $\phi  :V_1\to \RRR$ be a Morse function,
adjusted to $(\phi_1,u_1)$, and such that if
$p\in S(\phi)$
 then    $D_\delta(p)\subset\phi^{-1}(]a_k,a_{k+1}[)$
where
$k=\ind p$.                   Using deformations
along $(-u_1)$-trajectories it is easy to see that
$(V_1\setminus B_b^{(n-s-1)} ,A_b^{(s)})\sim
(\phi^{-1}(]a_0,a_{s+1}]), \phi^{-1}(]a_0,a_s])$.\quad$\square$
\enddemo

We proceed to construction of ranging pairs.
Let $v$ be an almost good
$f$-gradient.
Choose $\alpha,\beta\in]a,b[$ such that
$\alpha<\beta$
and that the intervals
$[a,\alpha], [\beta,b]$ are regular.
Let $\mu>0$ be so small that the map
$(x,\tau)\mapsto\gamma(x,\tau;v)$
is defined on $V_0\times [0,3\mu]$
and on
$V_1\times [-3\mu,0]$.
We denote these maps by
$\Psi_0:V_0\times [0,3\mu] @>>> W$, and, resp.
by
$\Psi_1:V_1\times [-3\mu,0]\to W$.
We assume further that $\mu$ is so small
that
$\Im\Psi_0\subset f^{-1}([a,\alpha[),~
\Im\Psi_1\subset f^{-1}(]\beta,b])$,
and that $\Psi_0$(resp.$\Psi_1$)-
induced metric
is 2-equivalent to the product metric.
 The $\Psi_i$-image of
$V_i\times [\alpha,\beta] (i=0,1)$ will be
denoted by $U_i(\alpha,\beta)$.

\proclaim{Lemma 4.12}
Let $\aaa ,\bb$ be $s$-submanifolds of $V_1$,~resp. $V_0$.
 Then there is an almost good $f$-gradient
$v_1$, $\smo$-close to $v$ and such that
\roster
\item $\supp (v-v_1)\subset \big(U_0(2\mu,3\mu)\cup U_1(-3\mu ,-2\mu)\big)$
\item $\kk _b(-v_1)\nmid \aaa$;
\qquad {\rm (3)} $\st v_1 (\aaa )\nmid \bb ~~,
  ~~ \stexp {(-v_1)} (\bb)\nmid \aaa$
\endroster
\endproclaim
\demo{Proof}
Note first that any $f$-gradient $v_1$, satisfying (1) is almost good.
Note also that the first property of (3) implies the second.
Now let $\xi\in \text{Vectt}(V_0)$ be a small vector field,
such that
$\Phi(\xi,\mu) (\kk _b(-v_1))\nmid \aaa$
(such vector field exists by 2.4). Then the vector field $v_0$
 which equals $v$ everywhere except $U_1(-3\mu,-2\mu)$
and in this
 neighborhood equals $(\Psi_0)_*((\xi)\times \frac {d}{dt})$,
satisfies (1) and (2).
Therefore the $s$-submanifold
$\st v_0 (\aaa )$ of $V_0$ is defined. Applying to $v_0$ the same
procedure as above (nearby $V_0$) we get an $f$-gradient
$v_1$, satisfying (1),(2) and the first part of (3).
$\qquad\square$
\enddemo

\proclaim{Theorem 4.13}
Let $\epsilon>0$. Then there is an almost good $f$-gradient
 $w$, and a ranging pair $(\VV _0,\VV _1)$ for $(f,w)$,
such that $\nr w-v\nr \leqslant\epsilon\aand
\supp (w-v)\subset U_0(\mu, 3\mu)\cup U_1(-3\mu,-\mu)$.
\endproclaim
\demo{Proof}
Choose $\epsilon'>0$ so small that $\epsilon'\leq \epsilon$
and $\epsilon'<\frac 1{\Vert \text{grad }f\Vert}\cdot
\min\limits_{x\in\Im\Psi_0\cup\Im\Psi_1}
df(v)(x)$.
Let $\VV _i =(\phi_i,u_i,\UU _i) (i=0,1)$ be
 an almost good $\MM$-flow on $V_i$, satisfying the conclusions
 of 1.19 with respect to $C=\epsilon'/4, \beta=\mu/3$.
Let $\delta_0>0$ be so small that $u_1,u_0$
are $\delta_0$-separated; then
$\delta_0$ is in the interval of definition
of $\KKK(-u_0), \KKK(u_1)$.
 Let $u$
be an almost good $f$-gradient, satisfying the conclusions of the
preceding lemma with respect to
 $\aaa =\kk (u_1)~,~\bb =\kk  (-u_0)$,
and such that $\Vert u-v\Vert<\epsilon'/2$.
Denote by $\widetilde W$ the cobordism $f^{-1}([\alpha,\beta])$
and by $\widetilde f, \widetilde u$
the restrictions of $f$, resp. $u$ to
$\widetilde W$.
Denote $f^{-1}(\alpha)$ by $\widetilde V_1$
and $f^{-1}(\beta)$ by $\widetilde V_0$.
Denote by $\widetilde{\AAA}$, resp.
$\widetilde{\BBB}$,
the image of $\KKK(u_1)$, resp. $\KKK(-u_0)$,
with respect to the diffeomorphism
$\stind ub\beta$, resp. $\stind {(-u)}a\alpha$.
Fix $]0,\delta_0[$ as the interval of defnition
of $\widetilde{\AAA}$, resp. $\widetilde{\BBB}$.
The cores of  $\widetilde{\AAA}$, resp. $\widetilde{\BBB}$
will be denoted by $\widetilde\aaa$, resp. $\widetilde\bb$.
Note  that
$\widetilde\aaa\nmid \stexp{(-\widetilde u)}(\widetilde \bb),~~
\widetilde\bb\nmid \stexp{(\widetilde u)}(\widetilde \aaa)$.
Choose $\delta_1>0$ such that $\widetilde u$ is
$\delta_1$ separated.
Choose $\theta>0$ so small that
 $\theta<\delta_0,~\theta<\delta_1,~\theta$
is in the interval of definition of
$\TTT(\widetilde\AAA,\widetilde u)$ and of
$\TTT(\widetilde\BBB,-\widetilde u)$ and
 for every $0\leq s\leq n-1$ we have:
$$
\TTT(\widetilde\AAA,\widetilde u)_{(\leq s)}(\theta)
\cap
\widetilde\BBB_{(\leq n-s-2)}(\theta)=\emptyset;\qquad
\TTT(\widetilde\BBB,-\widetilde u)_{(\leq s)}(\theta)
\cap
\widetilde\AAA_{(\leq n-s-2)}(\theta)=\emptyset\tag4.5
$$
Let $\WW _i=(\phi_i,u'_i,\UU'_i)$ (where $i=0,1$) be an
$\MM$-flow on $V_i$, subordinate to $\VV _i$ and satisfying
the condition 2) of 1.19 with respect to some
$\delta\in ]0,\theta[$. Let $\varkappa:[0,2\mu]\to[0,1]$
 be a $\smo$ function, such that $\supp\varkappa\subset ]\mu ,
2\mu[$ and $\varkappa(t)=1$ for
 $t\in [{\frac 43}\mu ,{\frac 53}\mu]$.

Define a vector field $w$ on $W$, setting
 $$
\gathered
w(x)=u(x) \ffor x\notin U_0(\mu,2\mu)\cup U_1(-2\mu,-\mu)\\
((\Psi_0^{-1})_* w)(y,t)=(\varkappa (t) u_0'(y),1)\ffor
y\in V_0,t\in [\mu,2\mu]\\
((\Psi_1^{-1})_* w)(y,t)=(\varkappa (-t) u_1'(y),1)\ffor
y\in V_1,t\in [-2\mu,-\mu].
\endgathered\tag4.6
$$
We claim that $w$ and $(\VV_0,\VV_1)$
 satisfy the conclusions of the theorem.
Our choice of $\epsilon'$
implies that $w$ is an $f$-gradient.
By the construction $w$ is almost good and
$\Vert w-v\Vert<\epsilon$.
The condition
$\supp (w-v)\subset U_0(\mu, 3\mu)\cup U_1(-3\mu,-\mu)$.
follows from (1) of 4.12 and from $(4.6)$ above.
To prove that
$(\VV_0,\VV_1)$
is a ranging pair for $(f,w)$ note first that
$0<\delta<\theta$,
 the gradients $u_0,u_1$ are $\theta$-separated.
Since $\widetilde u$ is
 $\theta$-separated,
$w$ is also
$\theta$-separated; (RP1) follows.
The formulas $(4.6)$ allow to give the following description
of the $(-w)$-trajectories, starting at a point $x\in V_1$.
\roster\item"(1)"
$\gamma(x,\cdot;-w)$ reaches $\widetilde V_1$ and intersects it
at a point
$x_1=\stind ub\beta \big(\gamma(x,\tau(x);-u_1')\big)$
where $\tau(x)\geq\mu/3$.
Then there are two possibilities:
\item"(2A)" $\gamma(x_1,\cdot;-\widetilde u)$
converges to some $p\in S(f)$.
Then the same is true for $\gamma(x,\cdot;-w)$.
\item"(2B)" $\gamma(x_1,\cdot;-\widetilde u)$
reaches $\widetilde V_0$ and intersects it at a point $x_2$.
Then the same is true for $\gamma(x,\cdot;-w)$
and we have:
\item"(3)"                        $
\gamma(x,\cdot;-w)$ reaches $V_0$ and intersects it at a point
$x_3=\gamma(\stind u\alpha{a} (x_2),\tau(x_2);-u_0')$,
where $\tau(x_2)\geq\mu/3$.            \endroster

(Similarly for $w$-trajectories starting at a point of $V_0$.)
This implies that $\stind wb\beta (\kk(u_1))=\widetilde\aaa$
and $\stind {(-w)} a\alpha (\kk(-u_0))=\widetilde\bb $, from
that one deduces easily that
$\kk_b(-w)\nmid \kk(u_1)$, and
$\kk_a(w)\nmid\kk(-u_0)$. Further, the
point (1) of the above description implies
that
$$
\forall\lambda<\delta_0\aand\forall s:0\leq s\leq n-1
\quad\text{we have}\quad
\stind wb\beta (\KKK(u_1)_{(s)}(\lambda))
\subset
\widetilde\AAA_{(s)}(\lambda)\tag*
$$
Therefore 2.17 applies to show that
$\theta$ is in the interval of definition
of $\TTT(\KKK(u_1),w)$.
Similarly, $\theta$ is in the interval of definition
of
$\TTT(\KKK(-u_0),-w)$.
Further, (*) and (4.5) imply that for every
$z\in \TTT(\KKK(u_1),           u)_{(\leq s)}(\theta)\cap \widetilde V_0$
we have:
$\stind u\alpha{a}
(z)\in V_0\setminus B_\theta(\indl {n-s-2},-u_0)$.
Therefore the third formula of $(4.6)$
together with (1.2) imply that
$\TTT(\KKK(u_1),w)_{(\leq s)}(\theta)
\cap V_0
\subset
B_\delta(\indl {n-s-2},-u_0)$
which imply (RP1) since $\delta<\theta$.
(RP2) is similar.\quad$\square$
\enddemo

\subhead{D. Equivariant ranging systems
and the proof of the main theorem}
\endsubhead
We assume here the terminology of Subsection B
of Introduction. It is easy to see that it suffices to prove
the Main Theorem in the case when
$[f]\in H^1(M,\ZZZ)$ is
indivisible, and
we make this assumption from now on.
This assumption implies that
$\forall z\in \bar M : F(z)-F(zt)=1$.
We assume also $S(f)\not=\emptyset$.
Choose a riemannian metric on $M$.
Then $\bar M$ obtains a $t$-invariant
riemannian metric.

\definition{Definition 4.14}
Let $u$ be any $f$-gradient.
Let $\Sigma$ be a non empty set of regular values
of $F$, satisfying $({\Cal S})$ below.

$({\Cal S})$:
 for every $A,B\in\RRR$ the set $\Sigma\cap\fpr F{[A,B]}$
 is finite;\quad
if
$\sigma\in\Sigma$ then
$\forall n\in\ZZZ :~~\sigma+n\in\Sigma$;\quad
if $\lambda,\mu\in\Sigma$ are adjacent, then there is only
 one critical value of $F$ between $\lambda$ and $\mu$.

 A set $\rans$ is called
\it
$t$-equivariant ranging system for $(F,u)$
\rm
, if
\roster
\item"(ERS1)" For every $\mu,\nu\in\Sigma,\mu<\nu$ we have:
$\{(A_\sigma , B_\sigma)\}_{\sigma\in\Sigma, \mu\leqslant\sigma\leqslant\nu}$
 is a ranging system for
$(F\vert\fpr F{[\mu,\nu]}, u)$.
\item"(ERS2)" $A_{\sigma -n} = A_\sigma\cdot t^n,\quad
B_{\sigma -n} = B_\sigma\cdot t^n$ ~for every $n\in \ZZZ$.
$\qquad\triangle$
\endroster
\enddefinition

Let $\rans$ be a $t$-equivariant ranging system for $(F,u)$.
 For $\nu,\mu\in\Sigma,~\nu<\mu$, denote
 by $H_{[\mu,\nu]}(u)$ the
homomorphism $H(u)$, associated (by  4.7) to the system
\break
$\{(A_\sigma , B_\sigma)\}_{\sigma\in\Sigma, \mu\leqslant\sigma \leqslant\nu}$.
For $\mu\in\Sigma$ set by definition
$H_{[\mu, \mu]}= \id :
H_*(F^{-1}(\mu)\setminus B_\mu, A_\mu)\to
H_*(F^{-1}(\mu)\setminus B_\mu, A_\mu)$.
 It follows from the construction that
$H_{[\mu-1,\nu-1]}(u)=
\mathbreak
t\circ H_{[\mu,\nu]}(u)\circ t^{-1}$ and that
$H_{[\nu,\theta]}(u)\circ H_{[\mu,\nu]}(u)=
H_{[\mu,\theta]}(u)$.
For $\nu\in\Sigma$ denote
\break
$t^{-1}\circ H_{[\nu,\nu-1]}(u)$ by $H_{(\nu)}(u)$.
It is an endomorphism of
$H_*(F^{-1}(\nu)\setminus B_\nu,A_\nu)$. We have obviously
$H_{
[\nu,\nu-k]}(u)
= t^k\circ\big(H_\nu (u)\big)^k$.
The next      lemma follows directly from      4.5 - 4.7.

\proclaim{Lemma 4.15}
Let $\rans$ be a $t$-equivariant ranging system
for $(F,u)$. Let $\nu,\mu\in\Sigma,\nu\leq\mu$ and
 $k\in\NNN$. Let $N$ be an oriented submanifold
of $\fpr F\mu\setminus B_\mu$ such
that $ N\setminus\text{\rm{ Int}}~ A_\mu$ is compact.
Let $L$ be a compact cooriented
submanifold of $\fpr F\nu\setminus A_\nu$.
Assume that $\dim N+\dim L=\dim M -1$.
Then:
\roster
\item $N'_k=\stind u\mu{\nu -k} (N)$ is an oriented
submanifold of $\fpr F{\nu -k}\setminus B_{\nu -k}$
 such that
$ N'_k\setminus\text{\rm{ Int}}~ A_{\nu-k}$
is compact.
If $N'_k\pitchfork L t^k$, then
$N'_k\cap L t^k$ is finite and
\newline
$N'_k ~\sharp ~ Lt^k=
\bigg\langle i^*( ]L[),\big(H_{(\nu)}(u)\big)^k\big(H_{[\mu,\nu]}(u)([N])\big)
\bigg\rangle$,
\newline
 where
$i:(F^{-1}(\nu)\setminus B_\nu , A_\nu )\to
(F^{-1}(\nu) ,F^{-1}(\nu)\setminus L)$ is the inclusion map.
\item For every $f$-gradient $w$, sufficiently close
 to $u$ in $C^0$-topology, $\rans$ is also
 a $t$-equivariant ranging system for $(F,w)$ and
$H_{(\nu)}(u)=H_{(\nu)}(w),
\break
H_{[\mu,\nu]}(u)= H_{[\mu,\nu]}(w). \qquad\square$
\endroster
\endproclaim
Passing to the proof of the Main Theorem
fix first two points $x,y\in S(f), \ind x=\ind y+1$,
and assume that
$F(\bar y)<F(\bar x)\leq F(\bar y)+1$.
Denote $\dim M$ by $n$; denote $\ind x$ by $l+1$,
then $\ind y=l$.
Choose some set $\Sigma$ of regular values
of $F$, satisfying $({\Cal S})$ of Definition 4.14.

Denote by $\theta$ the maximal element of $\Sigma$ with
$\theta< F(\bar x)$ and by $N(v)$ the
 intersection $D(\bar x,v)\cap F^{-1}(\theta)$;
$N(v)$ is an oriented submanifold of $F^{-1}(\theta)$,
diffeomorphic to $S^l$.
Denote by $\eta$ the minimal element of $\Sigma$,
satisfying $\eta>F(\bar y)$;
then $\eta\leq\theta<\eta+1$.
Denote by $L(-v)$ the
 intersection $D(\bar y,-v)\cap F^{-1}(\eta)$;
$L(-v)$ is a cooriented submanifold of $F^{-1}(\eta)$,
diffeomorphic to $S^{n-1-l}$.
Denote by $W$ the cobordism $F^{-1}([\eta,\eta+1])$.
Note that $\bar x\in \Wkr$. Denote
$F^{-1}(\eta)$ by $V_0$,
$F^{-1}(\eta+1)$ by $V_1$,
$\Sigma\cap[\eta,\eta+1]$ by $\Lambda$.

Denote by $\GG t_0(f;x,y)$ the subset of $\GG t(f)$,
consisting of all the $f$-gradients $v$, such that there is
an equivariant ranging system
$\rans$ for $(F,v)$ satisfying
$$\gather
N(v)\cap B_\theta =\emptyset,
\quad L(-v)\cap A_\eta=\emptyset, \tag4.7\\
H_*(F^{-1}(\eta)\setminus B_\eta,A_\eta)
\text{\rm~~ is a finitely generated abelian group.~~}\tag4.8\\
\endgather
$$
Now we shall prove 4 properties of the set
$\GG t_0(f;x,y)$.
\subsubhead{(1).
 $\GG t_0(f;x,y)$ is open in $C^0$ topology}
\endsubsubhead

This follows immediately from 4.15(2) and 1.10.
\subsubhead{(2).
 $\GG t_0(f;x,y)$ is dense in $C^0$ topology}
\endsubsubhead

By Th. 4.13 there is an almost good $F\vert W$-gradient
$w$ with $\Vert w-v\Vert<\epsilon/2$
and a ranging pair $(\VV_0,\VV_1)$ for
$(F\vert W,w\vert W)$.
By 4.11 this ranging pair generates the ranging
system
$\{(A^{(l)}_\lambda ,B^{(n-l-1)}_\lambda)\}_
{\lambda\in\Lambda}$.
Since $\supp (w-v)\subset \Wkr$, we can extend $w$ to a
$t$-invariant $F$-gradient on $\bar M$
(it will be denoted by the same letter $w$).
Checking through the proof of 4.13
shows that we can choose $\VV_1=\VV_0\cdot t$.
In this case for every $i$
 we have:
$A_\eta^{(i)}=A_{\eta+1}^{(i)}\cdot t$
and
$B_\eta^{(i)}=B_{\eta+1}^{(i)}\cdot t$;
therefore we can extend this ranging system
to a $t$-equivariant ranging system on $\bar M$.
The property (4.8) follows from
4.11 (3).
By definition
$N(w)=D(\bar x,w)\cap F^{-1}(\theta)$
is in
$A_\theta^{(l+1)}$,
and
$L(-w)=D(\bar y,-w)\cap F^{-1}(\eta)$
is in
$B_\eta^{(n-l)}$,
which (in view of 4.11(1))
implies (4.7) with respect to $w$.

Choose an $f$-gradient $\tilde w$,
with $\Vert w-\tilde w\Vert<\epsilon/2$
satisfying the transversality
assumption.
If only $\tilde w$ is close enough to
$w$, the system
$\{(A^{(l)}_\lambda ,B^{(n-l-1)}_\lambda)\}_
{\lambda\in\Lambda}$
is still a $t$-equivariant
ranging system for $\tilde w$ (by 4.15(2)), and
(4.7) still hold (by 1.10).

\subsubhead{3. For every $v\in\GG t_0(f;x,y)$ we have:
$n(x,y;v)=\frac {P(t)}{Q(t)}$ where
$P,Q\in\ZZZ[t]$
and $Q(0)=1$}
\endsubsubhead

Denote
$H_l(\fmin (\eta) \setminus B_{\eta},
A_\eta)$ by $\HH$. Denote
by  $[x]$ the element $[\stind v\theta\eta (N(v))]$
 of $\HH$. Denote by $[y]$ the element of
$\text{Hom} (\HH ,\ZZZ)$,
induced by the cohomology
class $i_v^*(]L(v)[)$, where
 $i_v: \big(F^{-1}(\eta )\setminus B_\eta ,
 A_\eta\big)
\hookrightarrow
\big(F^{-1}(\eta ) , F^{-1}(\eta )\setminus L(-v)\big)$
is
the inclusion map.
Denote by $h$ the endomorphism $H_{(\eta)}(v)$
of $\HH$.
By definition we have
$n_k(x,y;v)=
\big(D(\bar x,v)
\mathbreak
 \cap F^{-1}(\eta -k)\big)\sharp
\big(D(\bar y t^k, -v)\cap F^{-1}(\eta -k)\big)$. Therefore
4.15 implies             that
$n_k(x,y;v)=[y](h^k([x]))$ if $k\geq0$.
Since $F(\bar x)\leq F(\bar y)+1$,
we have $n_k(x,y;v)=0$
for $k<0$.
Therefore the demanded formula for
$n(x,y;v)$
follows immediately from
the next lemma.

\proclaim{Lemma 4.16}
Let $G$ be a finitely generated abelian group,
 $A$ be an endomorphism of $G$ and
$\lambda :G\to\ZZZ$
 be a homomorphism.
Then for every $p\in G$ the series
 $\sum_{k\geq 0} \lambda(A^k p)t^k\in \ZZZ[[t]]$
is a rational function of $t$ of
 the form $\frac {P(t)}{Q(t)}$,
where $P,Q$ are polynomials
and $Q(0)=1$.
\endproclaim
\demo{Proof }
It is sufficient to consider the
case of free f.g. abelian group $G$.
Consider  a free f.g.
$\ZZZ [[ t]]$-module
$R=G[[t]]$ and a homomorphism $\phi : R\to R$,
given by
$\phi = 1-At$. Then $\phi$ is invertible, the inverse homomorphism
given by the formula
$\phi^{-1}=\sum_{k\geq 0} A^k t^k$.
On the other hand the inverse of $\phi$
is given by Cramer
formulas,  which are rational functions with denominator
$Q(t)=\text{det}~(1-At). \qquad\square$
\enddemo

\subsubhead{(4).
 Let $v\in\GG t_0(f;x,y)$.
Let $U$ be a neighborhood of $S(f)$. Then there is $\epsilon>0$
such that for every $w\in\GG t_0(f;x,y)$
with
$\Vert w-v \Vert<\epsilon$
and $w(x)=v(x)$ for $x\in U$
we have:
$n(x,y;v)=n(x,y;w)$}\endsubsubhead

Let
$\rans$ be a $t$-equivariant ranging system
for $(F,v)$, satisfying (4.7) and (4.8).
Choose $\kappa<F(\bar x)$ so close to
$F(\bar x)$ that $D(\bar x,v)\cap F^{-1}([\kappa, F(\bar x)])
\subset U$. Denote $D(\bar x,v)\cap F^{-1}(\kappa)$
by $N_0$.
Choose $\epsilon'>0$
           so small that for
every
$w\in\text{\rm Vect }^1(W)$ with
$\Vert w-v\Vert<\epsilon'$
we have:
1) $w$ points
outward $F^{-1}([\theta,\kappa])$ in $F^{-1}(\kappa)$
and  inward $F^{-1}([\theta,\kappa])$ in $F^{-1}(\theta)$, and
2) every $(-w)$-trajectory starting at $N_0$
reaches $F^{-1}(\theta)$ and
$\stind w\kappa\theta
(N_0)\subset F^{-1}(\theta)\setminus B_\theta$
(see 5.6). Then the maps
$\stind w\kappa\theta \vert N_0,
\stind v\kappa\theta \vert N_0:
N_0\to F^{-1}(\theta)\setminus B_\theta$
are homotopic (via the homotopy
$(tw+(1-t)v)^{\rightsquigarrow}_{[\kappa,\theta]}$).
Therefore for every $f$-gradient $w$ with
$\Vert w-v\Vert<\epsilon'$ and $v=w$ in $U$ the elements
$[N(v)]=[\stind v\kappa\theta (N_0)]$
and
$[N(w)]=[\stind w\kappa\theta (N_0)]$
of $H_l(F^{-1}(\theta)\setminus B_\theta, A_\theta)$
are equal.
Similarly, there is $\epsilon''>0$ such that
for every $f$-gradient $w$ with
$\Vert w-v\Vert<\epsilon''$
and $w=v$ in $U$
the elements
$i^*_v(]L(-v)[)$
and
$i^*_w(]L(-w)[)$
of $H^*(F^{-1}(\eta)\setminus B_\eta, A_\eta)$
are equal. Now the proof  of
\it (4)\rm
is over by the
application of
4.15.

Note that another choice of liftings of $x,y$
to $\bar M$ leads to the incidence coefficient
$\tilde n(x,y;v)=n(x,y;v)t^k, k\in\ZZZ$.
Therefore
\it (4) \rm
  is true as it stands for another choice of liftings
and
\it (3) \rm
   turns to:
$\tilde n (x,y;v)=
\frac {P(t)}{Q(t)\cdot t^k}$,
where
$P,Q\in\ZZZ[t], Q(0)=1$.
Set $\GG t_0(f)$ to be the intersection
of  $\GG t_0(f;x,y)$ over all pairs
$x,y\in S(f)$ with $\ind x=\ind y+1$
and the proof of the Main Theorem
is over.\quad$\square$

\head{\S5. Appendix. $C^0$ perturbations
 of vector
 fields and their integral curves}\endhead

In this appendix we prove some technical
results on
integral curves of vector fields. The main
technical results
 are 5.1, 5.2, which state that the
trajectories of a $C^1$
vector field are in a sense stable
under small $C^0$
perturbations of the vector field.
For a manifold $M$ (without boundary)
 we denote by
$\hrrr$ (resp. $\vemm$) the vector space
of $C^1$ vector
fields on $M$ (resp. the vector space of
 $C^1$ vector
fields on $M$ with compact support).

\subhead{A. Manifolds without
boundary}\endsubhead
 In this subsection $M$ is a
riemannian manifold
without boundary, $v\in\vemm , n=\dim M$.

\proclaim{Proposition 5.1}
Let $a,b\in M, t_0\geq 0$
and $\gamma (a,t_0;v)=b$.

Then for every open
 neighborhood $U$ of
$\gamma(x,[0,t_0];v)$ and
every open
 neighborhood $R$ of $b$
there exist $\delta>0$ and an
open neighborhood
 $S\subset U$ of $a$ such that
$\forall x\in S$ and
$\forall w\in\vemm$ with
$\nr w-v\nr <\delta$ we have:
  $\gamma(x,t_0;w)\in R$ and
$\gamma(x,[0,t_0];w)\subset U$.
\endproclaim

\proclaim{Proposition 5.2}
Let $a,b\in M,t_0\geq 0$
and $\gamma (a,t_0;v)=b$.
Let $E$ be a submanifold
without boundary
of $M$
of codimension 1, such that
$b\in E$
and $v(b)\notin T_bE$.

Let $U$ be an open neighborhood
 of
$\gamma(a, [0,t_0];v)$ and $R$ be an
 open neighborhood of $b$
 Then for every $\theta>0$
sufficiently small
there exist $\delta>0$ and an open neighborhood
 $S\subset U$ of $a$, such that $\forall x\in S$
and $\forall w\in\vemm$ with $\nr w-v\nr<\delta$
we have:
\roster
\item $\gamma(x,[-\theta, t_0+\theta];w)\subset U$,
and $\gamma(x, [t_0-\theta,t_0+\theta];w)\subset R$.
\item There is a unique
$\tau_0=\tau_0(w,x)\in [t_0-\theta,t_0+\theta]$,
 such that
$\gamma(x,\tau_0;w)\in E$.
\item If $E$ is compact and $t_0$
is the unique $t$ from $[0,t_0]$
 such that $\gamma(x,t_0;v)\in E$,
 then $\forall y\in S$
the number $\tau_0(w,y)$ is the
unique $\tau$ from
$[-\theta,t_0+\theta]$ such that
$\gamma(y,\tau;w)\in E$.
\endroster
\endproclaim
\demo{Proof of 5.1}
The case $M=\RRR ^n$ (with the euclidean
metric)
 is obtained
immediately from
the next lemma.
For  $v\in \verr$
we denote by $\Vert v\Vert_1$
the norm of the derivative
$dv:\RRR^n\to L(\RRR^n ,\RRR^n)$.

\proclaim{Lemma 5.3}
Let $u,w\in\verr ,
\quad\nr u -w\nr<
\alpha,\quad\nr u\nr_1\leq D$,
where $D>0$.
Let $\gamma,\eta$
be trajectories of, respectively, $u,w$,
and assume
that $\vert \gamma(0)-\eta(0)\vert
\leq\epsilon$.
Then for every $t\geq 0$ we have:
 $\vert\gamma(t)-\eta(t)\mid\leq
\epsilon e^{Dt} +\frac \alpha D (e^{Dt} -1)$.
\endproclaim

\demo{Proof}
$\vert\eta'(t) -\gamma'(t)\vert =
\vert w(\eta(t))-u(\gamma(t))\vert
\leq
\vert w(\eta(t))-u(\eta(t))\vert +
\vert u(\eta(t))-u(\gamma(t))\vert
\leq
\alpha +D\cdot\vert \eta(t)-\gamma(t)\vert$.
Set $s(t)=\eta(t)-\gamma(t)$. Then
$\vert s'(t)\vert\leq\alpha +D\mods$,
and the standard argument
(see, for example, [2, p.117]),
shows that
$\vert s(t)\vert \leq\vert s(0)\vert\cdot e^{Dt}+
\frac \alpha D (e^{Dt} -1)$.
\quad
$\square$
\enddemo

Passing to the general case,
note that if $0\leq t_1\leq t_0$ and 5.1 is true
for the curve $\gamma (a,\cdot ;v)\vert [0,t_1]$
and for the curve
$\gamma (a',\cdot ;v)\vert [0, t_0-t_1]$,
where $a'=\gamma(a,t_0;v)$
then it is true for
$\gamma (a,\cdot ;v)\vert [0,t_0]$.
Therefore, applying successive subdivisions,
we can assume
that $\gamma (a, [0,t_0];v)$
belongs to the domain of a chart
$\phi :W\to\RRR^n$, and that $W\subset U$.
Choose  open neighborhoods
$W''\subset\overline{W''}\subset W'\subset
\overline{W'}\subset W$
of $\gamma (a,[0,t_0] ;v)$,
such that $\overline{W'}$ is compact,
and a $\smo$ function
$h:M\to [0,1]$ such that
$\supp h\subset W'$ and
$h(x)=1$ for $x\in W''$.
Choose $C>0$ so that
the metric, induced by $\phi$
from $M$ on $\RRR^n$ and the euclidean one
are $C$-equivalent in $\overline{W'}$.
For every $w\in\vemm$
denote the vector field
$\phi_*(h\cdot w)$ by
$\widetilde w \in\verr$.
Note that $\Vert\widetilde w\Vert_e
\leq C\Vert w\Vert_\rho$.
Let $R'\subset\phi (R\cap W'')$
be an open neighborhood of
$\phi (b)$.
By 5.3 there is $\delta_0>0$
and an open neighborhood $S'\subset \phi(W'')$ of
$\phi(a)$ such that for every
$u\in\verr$ with $\Vert \widetilde v - u\Vert_e<\delta_0$
and $x\in S'$ the
$u$-trajectory $\gamma (x,t;u)$
stays in $\phi(W'')$ for $t\in [0,t_0]$
and $\gamma (x,t_0;u)\subset R'$.
It is easy to see that the neighborhood
$\phi^{-1}(S')$ of $a$
and the number $\delta_0/C>0$
satisfy the conclusions of  5.1.~$\square$
\enddemo

\subsubhead{Proof of 5.2.
 1)~ The case $t_0=0$}
\endsubsubhead

We represent $\RRR^n$ as the product
$\RRR^1\times\RRR^{n-1}$;
the elements $z\in\RRR^n$ will be therefore
referred to as pairs $z=(x, y)$, where
$x\in\RRR , y \in \RRR^{n-1}$; $x$
is called  \it first coordinate of $z$. \rm
Set $\RRR_+^n=\{(x,y)\mid x>0\},~
\RRR_-^n=\{(x,y)\mid x<0\}$.
For $r>0$ denote
$]-r,r[\times B^{n-1}(0,r)$ by
$W_r$. Let $h:\RRR^n\to [0,1]$
be a $\smo$ function such that
$\supp h\subset W_2$ and $h(x)=1$
for $x\in W_1$.
 Choose a chart
$\phi:W\to \RRR^n$
of $M$, such that
$a\in W\subset U\cap R;~
\phi (a)=0;~
\phi (W\cap E)=\{0\}\times \RRR^{n-1}$
and
the first coordinate of
$\phi_*v (a)$ is $1$.
For a vector field
$w$ on $M$ we denote by $\widetilde w$
the vector field $h\cdot\phi_* w$ on
$\RRR^n$.
Choose  $C>0$ so that
the metric, induced by $\phi$
from $M$ on $\RRR^n$ and the euclidean one
are $C$-equivalent in $\overline{W_2}$.
Then for every $w\in\vemm$ we have
$\Vert\widetilde w \Vert_e\leq C\Vert w\Vert_\rho$.
Let $r\in]0,1[$ be so small that the first coordinate
of $\phi_*v(x)$ is not less than
$1/2$ for $x\in W_r$.
Let $\theta>0$ be so small
that
$$
\gamma (0, [-\theta ,\theta];\widetilde v)
\subset W_r ,\quad
\gamma (0, \theta;\widetilde v)
\in     \RRR^n_+
      \aand
\gamma (0, -\theta ;\widetilde v)
\in \RRR^n_-.
$$
By 5.3 there is $\delta_0>0$  and a neighborhood
$S_0\subset W_r$ such that
$$
\gamma (x, [-\theta ,\theta]; u)
\subset W_r ,\quad
\gamma (x, \theta; u)
\in     \RRR^n_+
      \aand
\gamma (x, -\theta ; u)
\in \RRR^n_-
$$
whenever $x\in S_0,~ u\in\verr,
\Vert u-\widetilde v\Vert_e<\delta_0$.
Therefore for $x\in S_0,~ u\in\verr,
\Vert u-\widetilde v\Vert_e<\delta_0$
there is $\tau_0=\tau_0(x,u)\in ]-\theta,\theta[$
such that $\gamma(x,\tau_0)\in
\{0\}\times \RRR^{n-1}$.
Further, if $\Vert u-\widetilde v\Vert_e<\delta_1
=\min (\delta_0,1/2)$
the first coordinate of $u$ is positive
in $W_r$,
and therefore for every
$x\in S_0$, there is one and only one
$\tau_0=\tau_0(x,u)\in[-\theta,\theta]$
with $\gamma(x,\tau_0;u)\in\{0\}\times\RRR^{n-1}$.
It is easy to see that
(1) and (2) of the conclusions of our
proposition
hold for $\delta =\delta_1/C$
and $S=\phi^{-1}(S_0)$.
 (3) is an immediate
consequence of (2).

\subsubhead{2) ~ General case}
\endsubsubhead

The part 1) of the present
proof  applied to the point $b$
and $t_0=0$, implies that there is
$\theta_0>0$ such that for every $\theta\in]0,\theta_0[$
there is $\delta_0(\theta)>0$ and a neighborhood
$R_0(\theta)\subset U\cap R$
such that for every $x\in R_0(\theta)$ and for every
$w\in\vemm$ with
$\Vert w-v\Vert<\delta_0(\theta)$
we have:
$$\gather
\gamma (x,~[-\theta,\theta];w)\subset U\tag5.1\\
\exists ~!~
\tau_0=\tau_0(w,x)\in[-\theta,\theta],
\text{\rm such that }
\gamma(x,~\tau_0;w)\in E\tag5.2
\endgather
$$
Choose $\theta_1>0$ such that
$\gamma(x,[-\theta_1, t_0+\theta_1];v)\subset U$,
and that $\gamma(x,[t_0-\theta_1,t_0+\theta_1];w)\subset R$.
Applying 5.1 one obtains easily
that
for $\theta\in]0,\theta_0[$
there exist
$\delta_1(\theta)>0$
and an open neighborhood $S(\theta)\subset U$ of $a$,
such that for every $y\in S(\theta)$ and every
$w\in\vemm$ with $\Vert v-w\Vert<\delta_1(\theta)$ we have
$$\gather
\gamma(y,[-\theta_1,~t_0+\theta_1];w)\subset U
\aand \gamma(y,[t_0-\theta_1,t_0+\theta_1];w)\subset R
\tag 5.3\\
\gamma(y,t_0;w)\in     R_0(\theta)\tag 5.4
\endgather
$$

We claim that for every $\theta\in]0,\min (\theta_0,\theta_1)[$
the number $\delta_2(\theta)=
\min (\delta_1(\theta),\delta_0(\theta))$
and the neighborhood $S(\theta)$ of $a$ satisfy
(1) and (2) of conclusions of 5.2.
Indeed, (1) is immediate from (5.3).
To prove (2), note that the existence of the unique
$\tau_0\in [t_0-\theta,t_0+\theta]$ such that
$\gamma(y,\tau_0;w)\in E$
is equivalent to the existence of the unique
$\tau_1\in [-\theta,\theta]$ such that
$\gamma( \gamma(y,t_0;w),\tau_1;w)
\in E$
which is guaranteed for $y\in S(\theta)$ and
$\Vert w-v\Vert\leq\delta_2(\theta)$
by (5.4) and (5.2). Therefore for non compact
$E$ the proof is over.
For the case of compact $E$ choose $\theta_3>0$ such that
$\gamma(a,[-\theta_3,t_0[,v)\cap E=\emptyset$.
For every $\theta\in]0,\theta_3[$ choose
$\delta_3(\theta)>0$ and a neighborhood
$S_3(\theta)$ of $a$, such that
$\gamma(x,[-\theta_3,t_0-\theta];w)\subset M\setminus E$
for $x\in S_3(\theta)$ and $w\in\vemm$ with
$\Vert w-v\Vert <\delta_3(\theta)$.
Then for every $\theta\in ]0,\min (\theta_1, \theta_2, \theta_3)[$
the number $\min(\delta_1(\theta),
\delta_2(\theta),
\delta_3(\theta))$,
and the neighborhood
$S_3(\theta)\cap S(\theta)$ of $a$
satisfy (1) - (3) of 5.2.
\quad $\square$

\subhead
{B. Manifolds with boundary}
\endsubhead
Let $W$ be a compact riemannian
 manifold
with boundary. Recall from [19, \S 8C]
that we denote by
 $\text{\rm Vect}^1 (W,\bot)$
 the space of $C^1$ vector fields $v$ on $W$,
such that $v(x)\notin T_x(\partial W)$
for $x\in\partial W$.
         Choose an embedding of $W$ into a closed
manifold
$M$ without boundary, $\dim M =\dim W$
(for example,one can take the double of $W$).
Pick a riemannian metric on $M$
extending that of $W$ (the existence
of such a metric is easily proved using
the standard partition of unity argument).
The same argument proves that every $C^1$
 vector
field on $W$           can be extended
to a vector field $\widetilde v\in \text{\rm Vect }^1(M)$
such that $\Vert \widetilde v\Vert\leq 2\Vert v\Vert$.
 The following corollaries are "$C^0$-analogs"
of Propositions 8.10, 8.11 of [19].
They are proved similarly to 8.10, 8.11 of [19]
using Propositions 5.1, 5.2 of the present appendix
instead of  [19, Prop. 8.2].
\proclaim{Corollary 5.4}
Let $v\in\vew , x\in W$. Assume
that the
$v$-trajectory $\gamma(x,\cdot~;v)$
reaches
the boundary at a moment $T\geq 0$.
 Let $U$
be an open neighborhood of
 $\gamma(x,[0,T];v)$,
and $R$ be an open neighborhood
of $\gamma(x,T;v)$.
Then there is $\delta>0$ and
a neighborhood $S$ of
$x$, such that for every $y\in S$
 and for every
$ w\in\vew$ with
 $\nr w-v\nr<\delta$ the
 trajectory
$\gamma(y,\cdot ;w)$
reaches the boundary at a
moment $T(y,w)\geq 0$, and
$\gamma(y,[0,T(y,w)];w)
\subset U$
and $\gamma(y,T(y,w);w)\subset R.
\quad\square$
\endproclaim

\proclaim{Corollary 5.5}
Let $v\in\vew ,x\in W, ~T\geq 0$. Assume,
that $\gamma (x,T;v)
\in\overset\circ\to W$.
 Let $U$  be an open
neighborhood of
$\gamma(x,[0,T];v)$, and $R$
be an open neighborhood of
 $\gamma(x,T;v)$.
Then there is $\delta>0$
and a neighborhood
$S$ of $x$, such that for every $ y\in S$
 and for every
$ w\in\vew $ with $\nr w-v\nr<\delta$
the trajectory
 $\gamma(y,\cdot ;w)$ is defined
on $[0,T]$
and $\gamma(y,[0,T];w)
\subset U$ and
$\gamma(y,T;w)\in R$.
\quad$\square$
\endproclaim

Let $v\in\vew$. Denote by
$V_0$, resp. by $V_1$ the set
of $x\in\partial W$, where
$v(x)$ points inward $W$, resp.
outward $W$. For a subset $Z\subset W$
denote by $\tau(Z,v)$ the set
$\{\gamma(z,t;v)\vert z\in Z, t\geq 0\}$.
The next corollary is deduced from 5.4
by an easy compactness argument.

\proclaim{Corollary 5.6}
Let $K\subset W$ be a compact such
that every $v$-trajectory starting at a point of $K$
reaches the boundary.
Let $U$ be an open
neighborhood of $\tau(K,v)$.
Let $R\subset V_1$ be an open
neighborhood of $\tau(K,v)\cap V_1$.
Then there is $\delta>0$ such that
for every $w\in\vew$ with $\Vert w-v\Vert<\delta$
each $w$-trajectory starting
at a point of $K$ reaches the boundary and
we have:
$\tau(K,w)\subset U,\quad \tau(K,w)\cap V_1\subset R$.
\quad$\square$
\endproclaim

\enddocument